\documentclass[12pt,reqno]{amsart}
\usepackage[all]{xy}
\usepackage[margin=0.9in]{geometry} 
\usepackage{amssymb,amsmath,amsthm}
\numberwithin{equation}{section}
\usepackage{fancyhdr}
\newcommand{\R}{\mathbb{R}}

\theoremstyle{definition}

\newcommand{\ben}{\begin{enumerate}}
\newcommand{\een}{\end{enumerate}}
\newcommand{\eit}{\begin{itemize}}
\newcommand{\beq}{\begin{equation}}
\newcommand{\eeq}{\end{equation}}

\renewcommand{\leq}{\leqslant}
\renewcommand{\geq}{\geqslant}

\begin{document}

\title{Accepted proofs:  Objective truth, or  culturally robust?}
\author{Andrew Granville}
\dedicatory{\, \hfill Vanity of vanities! Vanity of mathematics! \\ \, \hfill -- Frederick II of Prussia (1778)}

 \address{D{\'e}partment  de Math{\'e}matiques et Statistique,   Universit{\'e} de Montr{\'e}al, CP 6128 succ Centre-Ville, Montr{\'e}al, QC  H3C 3J7, Canada.}
   \email{andrew.granville@umontreal.ca}

 \thanks{The quotation in the dedication is extracted from a letter \cite{Eck} from Frederick II to Voltaire, complaining about the ``expert advice'' given to him by one Leonhard Euler. A special thanks to Michael Hallett for explanations of the history of the development of axiomatic systems, help with formulating my viewpoints and pointing out some remarkable flows of ideas, Rodrigo Ochigame for developing similar ideas from a very different and thought-provoking perspective, as well as Peter Scholze for several comments and corrections.
 I would also like to thank Dave Benson, Kevin Buzzard, Dan Goldston, Deirdre Haskell, Anthony Kennedy, Heather Macbeth, J\'an Min\'a\v c, Mel Nathanson,  Peter Sarnak,  and  Carlos Simpson for helpful discussions, responding to emails and useful references, as well as the referee for a careful reading and good suggestions.
 The author is partially supported by grants from NSERC (Canada). }  
\begin{abstract}
How does the mathematical community  accept that a given proof is correct? 
Is  objective verification based on explicit axioms feasible, or must the reviewer's experiences and prejudices necessarily come into play?
Can automated provers avoid  mistakes (as well as experiences and prejudices) to provide objective verification? And can an   automated prover's claims be provably verified?

We will follow examples of proofs that were found to be flawed, but then corrected (as the proof plan was sufficiently robust), as well as accepted ``proofs'' that  turned out to be fundamentally wrong. 
What does this imply about the desirability of the current community standard for proofs?

We will discuss whether mathematical culture is unavoidably part of the acceptance of a proof, no matter how much we try to develop foolproof, objective  ``proof systems''.  
\end{abstract}

\maketitle
 
 
\vskip-.2in

\centerline{ \textsc{Contents}}
{\obeylines
\noindent 1. \emph{Proof -- why and how}:  Hilbert's axiomatic approach to mathematical theory.
\noindent 2.  \emph{Living with, and ignoring, the G\"odel crisis}
\noindent 3.  \emph{Formal proof vs  culturally appropriate, intuitive explanation.} 
Including a discussion of  approaches to proof and the use of language.
\noindent 4.  \emph{What is an accepted proof in pure mathematics?}   The realities of the refereeing process, and   what is  ``interesting''.
\noindent 5.  \emph{Mistakes}. Famous recent  mistakes   and how the community responded to them.
\noindent 6.  \emph{Rethinking axioms and  language}.  Investing in new theories has its pitfalls.
\noindent 7.  \emph{Computers and proofs}.  Establishing, assisting with, and generating, proofs.
\noindent 8.  \emph{Uses of computers in major theorems}. Inspiring examples and  computers   in their proofs.
\noindent 9.  \emph{Computer error}. The obvious pitfalls and whether we can elude them.
\noindent 10.  \emph{Protocols for automated theorem checkers/provers}. What can we hope for?
\noindent 11.  \emph{The future of proof}.  Computer proofs and a modified Turing test.
\noindent 12.  \emph{The Lean Theorem prover}. Extraordinary recent developments with  Peter Scholze.
\noindent 13.  \emph{Myths of objectivity}.  Growing myths about the correctness of computer-proofs fit with well-known examples of refuted objectivity.
\noindent 14.  \emph{Will machines change  accepted  proof}? Some conservative speculation about the future.

}

\section{Proof -- why and how}

To begin we discuss why proofs are desirable, what is the generally accepted approach to proof, and what aspects are theoretically problematic.
\subsection*{The purpose of proof}
Aristotle wrote
\begin{quote}
If ... understanding is as we posited, it is necessary for demonstrative understanding ... to depend on things which are true and primitive and immediate and more familiar than and prior to and explanatory of the conclusion.
\end{quote}
More generally we can prove new concepts by reducing them to those which have already been accepted;
in particular one does not need to always deduce the latest assertion directly from the axioms, as beautifully  explained  by Nathanson \cite{Nath08}:
\begin{quote}
We mathematicians have a naive belief in truth. We prove theorems. Theorems are deductions from axioms. Each line in a proof is a simple consequence of the previous lines of the proof, or of previously proved theorems.  \emph{Our conclusions are true, unconditionally and eternally.}  
\end{quote}
That is, modern mathematicians   dream that all theorems should be provable from appropriate axioms, and  we just have to find the proofs.
However, even back in the late 19th and early 20th century, it was found to be difficult to decide what precisely is meant by these statements.

\subsection*{Hilbert's vision} In early discussion of the foundations of mathematics, the consensus was to  build on the already-acknowledged correct 
``Aristotlian primitives'', rather than question whether one could (and should) start with any set of coherent, consistent axioms. For example, Hallett \cite{Hall} writes:
\begin{quote}
Like Dedekind, Cantor argued that progress in mathematics depends on conceptual innovation, the central constraints being the `integration' of the new concepts with already accepted concepts and the condition that the concepts be consistent. 
\end{quote}
Although this suggests that there is unlikely to ever be an eternally fixed set of axioms, there  were two important attempts to formulate a single sufficient axiomatic system: the first, due to Frege,  was flawed because it was inconsistent, giving rise to  Bertrand Russell's \emph{antinomies}, that is, pairs of reasonable-sounding statements which are contradictory (leading to logical paradoxes\footnote{For example, the wonderfully relevant paradoxical assertion, ``This statement is not true''.}; though
G\"odel \cite{God} points out that these antimonies ``do not appear  within  mathematics'' itself but rather within its ``interpretation'')
The second, developed in  Russell and Whitehead's \emph{Principia Mathematica}, was also flawed since it contained three axioms that they could not justify (and indeed seemed ad hoc).

Cantor's work on infinities was partly developed under the assumption\footnote{Which Cantor called  ``a law of thought'' in 1883, what might be called ``self-evident'' today.} that every set can be well-ordered.  Better understanding this assumption formed part of Hilbert's First Problem, resolved in 1904 by Zermelo's   \emph{Axiom of choice}. These issues led Hilbert in the 1920s to develop his foundational program, for the first time explicitly beginning with axioms to establish a framework to study the foundations of mathematics, where abstract, ideal mathematics forms the very subject matter, without earlier prejudices that a particular route is desirable.  The axioms create the theory, not vice-versa.

Hilbert wanted axioms that are developed consistently through clear rules, and do not lead to antinomies and paradoxes.\footnote{Hallett \cite{Hall} remarks that part of the inspiration for this came from the proofs of the consistency of the alternative, non-Euclidean, geometry of Gauss, Bolyai and Lobatchevsky.} Hilbert therefore did not pretend, like Frege and  Russell, to select  axioms   because they are `truths', but rather because they are reasonably straightforward assertions which allow the subject matter to develop consistently.  In particular an axiomatic system does not need to start with   `basic incontrovertible truths', but rather  starting points that make sense and are consistent.
\begin{quote}
Hilbert did not share Russell's interest  in `getting the foundations right', in isolating the `right' set of primitives. For Hilbert, there is no `right' set of primitives; some might be better than others for certain purposes, but there is never really any final word, even for relatively simple theories \footnote{Other leading logicians, like Frege, felt that a line is a line, the physical entity in the common vernacular, and if you produce a theory yielding something different with an obscure, though consistent set of axioms, then you should discard that  theory; Hilbert disagreed.} ... Hilbert stressed that any theory is only a `schema of concepts'; it is `up to us' how to fill it with content.
\hfill Hallett \cite{Hall}
\end{quote}
Indeed one prefers to supplement the basic  theory, not so much by using new axioms, but rather through  different conceptual frameworks, embracing a much wider range of concepts than just the numbers and their properties. Axioms cannot  be proved,  and are not in need of proof.

Hilbert also worried about language and interpretation, demanding that all interpretations  of the theory should be isomorphic, and all deductive techniques should be invariant under different interpretations.\footnote{Not only should all models of the same theory be isomorphic, but no matter how one expresses an idea, and no matter how one reads it, the ideas must all boil down to the same thing.}  Hilbert's program lies at the base of modern mathematics so let's clarify the details of some of  the ingredients in these theories.

\subsection*{Constructing a ``formal reasoning system''}  Language is often imprecise, and people's interpretations and recollections can differ.  We want  to create a language that promotes accuracy and a lack of ambiguity, and so we focus on its rules:

We begin with an alphabet of symbols and variables; anything you like   as long as it is countable and then sentences are only allowed to be finite combinations of those symbols.
Then we have grammatical rules of how to write meaningful sentences; again these can be whatever you like as long as they seem to be consistent and flexible.  The only substantial idea is that there should be  ``theorems'',   formulas without free variables (also called ``closed formulas'').  The first theorems are called axioms and can be enunciated without proof; then all others must be deduced using (only) our deduction rules from the axioms. Moreover one should be able to decide in a finite time (defined appropriately) whether a proof is correct.
 Consistency is key: we should not be able to determine a theorem and its negation from our axioms.  
 
 \emph{First order logic} keeps the special symbols down to a sensible minimum, allowing the quantifiers $\forall, \exists$, the formula connectives $\lor, \land$ and $\implies$ (with which we can express $A$ implies $B$ by $A\implies B$), and the negation $\neg$. This is enough of a  language to  work with most intuitive mathematics, and the hope was that it, or it supplemented by one or two necessary refinements, could cover all mathematical truths.\footnote{Codifying the means allowed in proof has led to modern logic, based on a small list of logical primitives, and such simple rules of inference.}   But we must always remember that
 \begin{quote}
Mathematics is writing. For all the quantification it makes possible and all the technological and scientific discoveries it has helped to produce, it is ultimately words upon words. There is a bedrock of definitions ... crosscut by axioms ... whose only restrictions are that they not contradict one another. From this starting material we derive the terse assertions of consequences that are known as theorems and lemmas and corollaries. ... The arguments are the key. These are the ``proofs'' ... without which the assertions are just so much blather. The proofs actually conjure mathematics into existence.\newline
. \hfill -- Dan Rockmore  \cite{Rock}, 2022.
 \end{quote}
\subsection*{Plurality} 
 
The foundations  of mathematics starts with a set of axioms, but which set?  Hilbert's proposals leave open the possibility of  working within different axiomatic systems, and perhaps those different axiomatic systems will lead to contradictory conclusions to the same simple questions. Then how do we decide which axiomatic system is the correct one to use?
This has not really been an issue for research mathematicians who  accept that one cannot avoid  \emph{deductive pluralism} (see  Hosack \cite{Hosa}), as long as each axiomatic system proceeds consistently. Let's recall two famous examples:

Zermelo introduced the Axiom of Choice (AC) in 1904, though it was only in the late 1930s that G\"odel was able to show it is consistent with the other axioms.\footnote{For the experts, by the ``other axioms", I mean  Zermelo-Fraenkel set theory.}  However the shock came in 1963 when Paul Cohen showed that it is fully independent from the other axioms since $\neg$AC is also consistent with the other axioms. Therefore do we work with the axiom of choice, or not? 
 
 Then there is Cantor's 1878 \emph{Continuum Hypothesis}, the claim that there are no infinities that in size\footnote{Size here is defined in terms of 1-to-1 correspondences.} lie strictly between the set of integers and the set of  reals. In 1940 G\"odel showed that one could not use the other axioms (including the Axiom of Choice) to show that there are such infinities, and in 1963 Cohen showed that 
the other axioms could not show that there are no such infinities. That is, the Continuum Hypothesis is independent of the standard axioms.\footnote{Thus if ZFC denotes Zermelo--Fraenkel set theory together with the axiom of choice then 
\emph{both}  ZFC+CH and ZFC+$\neg$CH are  consistent.}
 
 Research mathematicians have little difficulty in accepting this plurality.  Indeed some of the  most interesting research programs in recent years in pure mathematics have needed to rethink some of the starting assumptions and definitions to make further progress on questions of interest (as we will discuss later). 
   Nathanson \cite{Nath09} remarks that this plurality can lead to the misguided notion that 
  \begin{quote}
  mathematics is the logical game of deducing conclusions from interesting but arbitrarily chosen finite sets of axioms and rules of inference. 
  \end{quote}
  Indeed even G\"odel \cite{God} agreed that 
  \begin{quote}The truth of the axioms from which mathematics starts out cannot be justified, and therefore [one can argue] the drawing of consequences
... is a game with symbols according to certain rules, not insights.
  \end{quote}
  Nonetheless there is broad agreement to begin with the nine axioms of ZF, perhaps adding the axiom of choice.
  To be widely accepted, any radically new starting points need to be selected for good reasons, really important reasons, or they will gain little traction. 
  For example,  there is much research today on which axiom or axioms need to be added to the standard ones to ``decide'' the truth of the Continuum Hypothesis.\footnote{Woodin's (*)-axiom, or Martin's maximum for forcing, or variants of either, and whatever we choose needs to lead to a rich and consistent theory of infinities which can answer many of the important questions.}
  Thus although Hilbert's dream has been widely accepted,  G\"odel \cite{God} noted that we still make choices as argued by Dedekind and Cantor, just at a different stage and with a different level of confidence that what we do is coherent and consistent.

\subsection*{A key goal} 
Axiomatic systems are not only there to ensure correctness, but hopefully to provide a framework to prove all the theorems that are worth proving. After all if we can state a theorem, surely it should be provable in our ``theory'', in that the theory should not have any artificial limitations.
  \begin{quote}
 Hilbert ... demand[ed] that an axiomatised theory be complete, meaning that the axiom system be able to derive all the important facts, or all facts of a certain sort.
 \hfill Hallett \cite{Hall}
\end{quote}
As we wrote above, we wish to create a (finitely described) language/theory to quantify and work with all intuitive mathematics, perhaps enhanced by further axioms as we determine new issues that are   independent of the axioms we are already using, so that any true theorem is accessible and provable.   Hilbert's dream was  that one could create axioms and a theory that would allow one to prove or disprove any given mathematical claim.

In 1931 G\"odel showed that this dream is inescapably impossible.

\section{Living with, and ignoring, the G\"odel crisis} 

In 1931 Kurt G\"odel shook mathematics (literally) to its foundations with his \emph{incompleteness theorems}:

I) No consistent finite set of axioms and rules\footnote{I am avoiding necessary subtleties. To be more precise, G\"odel's theorem can be stated as: The set of provable statements, although not recursive, is recursively enumerable, whereas the set of true statements is not.} can be used to prove all true theorems about the integers.

\noindent Even worse:

II) No consistent finite set of axioms and rules can prove itself to be consistent!
\smallskip

This absolutely contradicts Hilbert's dream that one could axiomatize so as to prove or disprove any given mathematical claim.
The first is frustrating, while the second unveils inescapable limitations in the possible formulation of axiomatic theories and indeed mathematical foundations:
\begin{quote}
Formal, axiomatic set theory ... cannot be a final foundational theory for mathematics, [yet] nothing else [can] be.
 \hfill Hallett \cite{Hall}
\end{quote}
Or as von Neumann put it in 1930:
\begin{quote} 
There can be no rigorous justification for classical mathematics.
\end{quote}

So how do mathematicians deal with this existential crisis in their subject?  The only answer is that they learn to live with it.\footnote{From Janella Baxter on Twitter (@BaxterJanella, Oct 6, 2021):``When my husband was a math PhD [student], G\"odel's theorem terrified him. He worried the dissertation topic he'd selected was impossible to solve. He made it through, but he's baffled why more mathematicians don't share this anxiety. I think there's a philosophy of math paper to be written...''} Here's how:

\subsection*{How to (not) deal with incompleteness}  \mbox{}

At the beginning there are axioms;  consistent, independent and all powerful. 

They are designed to be used for the rest of eternity to prove theorems, building up from these axioms a towering, self-justified, self-supporting edifice, a structure that allows one to objectively  create a theory of everything.

However humans have found that there are questions that those axioms could not answer, so   had little option but to supplement the axioms
with further axioms to answer those previously unanswerable questions. But,  no matter how many axioms were added it seemed that not all questions could be answered, not all theorems could be proved, some theorem always seemed to be undecidable using just the given axioms. 

Indeed it came to pass that G\"odel established  there would always be mysteries, no matter how many or what axioms or language one starts with,  questions that cannot be answered, created out of the very  axioms that one starts with, in the language one has chosen.

\emph{How should one handle  this impasse?}  What are theorem tower  builders     to do? 

When one occupies the penthouse suite of a very high and seemingly solidly built tower, and there is a small fire below but it is difficult to review the damage done,   one has   options:

$\bullet$\ One can abandon the tower and try to rebuild it; if one can only work out how to do so given the new fire code. This is risky as it may mean  never ascending as high again.\footnote{And who knows? They might change the fire code halfway through reconstruction, adding delays and extra cost.} 

$\bullet$\   After  the smoke and noise have mostly dissipated one can  continue on as if nothing is amiss.  Maybe the cinders will smoulder for a bit and then go out of their own accord? And hey, if the building collapses, which seems unlikely given how well built it mostly is, then we won't be around to learn more.

The latter option is the choice of most pure mathematicians to cope with the G\"odel crisis.
Investing in the building that is, not in the building that may be yet to come.


\subsection*{Inescapable logic}  For many mathematicians, G\"odel's objections are melodic whistling in the wind -- perhaps practical questions in algebra and analysis  would never face these sorts of epistemological problems? Surely we need  a formal reasoning  system that is reliable in all ``reasonable circumstances'' despite being incomplete (though seemingly incomplete only when one is truly looking for confounding problems).  Mathematicians floated on that cloud until the devastating Ramsey theorem of Paris and Harrington \cite{PaHa77}: Although only mildly more complicated than what people work with every day in Ramsey theory,   it cannot be proved in Peano arithmetic since, roughly speaking, it implies that Peano arithmetic is consistent (and  G\"odel's  second incompleteness theorem says that you can't use a theorem in Peano arithmetic to prove its consistency).  There is now a collection of beautiful and natural theorems with the same property.
 
 The general mathematical culture is to not worry about these things too much. If one works on a question close   to where some  problem has been found that is undecidable within the ``standard'' axiomatic framework then one needs to act carefully,\footnote{By ``undecidable" we mean that if it is true it is not provable using the axioms in the theory, and if it is false it is not refutable using those same axioms.}  but for the most part it seems like a distant problem, and not one that we meet in our day-to-day work.\footnote{One does hear people suggest that popular unsolved problems might be ``undecidable'' within our axiomatic framework, which I regard as hubris; just because we have not yet found a good understanding of something does not make it an eternal mystery.  Knuth \cite{Knu} even makes the tenuous argument ``the Goldbach conjecture ... [is] a problem that's never going to be solved. I think it might not even have a proof. It might be one of the unprovable theorems that G\"odel showed exist ... we now know that in some sense almost all correct statements about mathematics are unprovable,'' and goes on to claim that Goldbach must be ``true because it can't be false'' for which he then gives a standard heuristic.  The only salvageable truth from this is that Goldbach might not be provable in Peano arithmetic since there might be a different model of integers that satisfy Peano arithmetic yet for which Goldbach fails; however if so then we tinker with our axioms and add one to ensure we remain in the usual integers and then  Goldbach should be provable. People have  made analogous fatuous claims about the Riemann Hypothesis, the twin prime conjecture, etc with no real substance to back their claims. There is a good discussion of all this at  \centerline{https://mathoverflow.net/questions/27755/}. }
 As Peter Sarnak likes to say of mathematicians, ``we are working people''.\footnote{My own attitude to this is akin to my attitude to food sourcing: I would like to eat wholesome, healthy, unadulterated, non-exploitative food, and I make choices at the supermarket to do so by taking a superficial look at the labels. However I am quite aware that some of their choices, and their labelling, may be questionable. Indeed I do not really believe my local organic supermarket to be ``trusted authorities'', but I do not want to spend my time carefully following up each purchase, holding my purveyors accountable,  unless something in particular makes me suspicious (or naseous).}
Indeed for the last century the basic axiomatic system  (ZFC) has remained accepted as the essentially unchanged foundation for most of modern pure mathematics,   surprising in the light of G\"odel's results, and yet it works.

That there is no way to evade G\"odel's incompleteness theorems is so counterintuitive to the prevailing culture that
 some top people (including two Fields medallists) have suggested that somehow one can evade these issues of what is provable by writing suitably capable computer programs.  Voevodsk\u{\i}'s evasion is that it would be enough to have a ``reliable'' system, but is  vague on what he means  \cite{Voev}.\footnote{Presumably  the idea is to define a plausible concept of ``formal reasoning system'' that allows one to evade G\"odel's result.} Thurston \cite{Th84} seems to suggest that it is difficult to hide errors and have a program run since one has to be syntactically correct, which seems to confuse two different issues: most mathematical papers that contain important mistakes are written in a coherent manner, obeying the syntax and grammar of our system of communication -- a misinterpretation of one of these is rarely where the mistake lies (though see the Biss case discussed below).  
In 1950 Turing \cite{Tur1} had already remarked that  one can adapt G\"odel's proof to create questions that any given computer language cannot answer, in that the question is undecidable within the computer language's design.  
 
 One can interpret the \emph{Church-Turing thesis} as the working hypothesis that all sensible computational systems are ``equivalent'' and ``universal''. That is, they can each calculate anything that is calculable, though time may be an issue. Thus there is  theoretically little difference between the capability of different programming languages (although they might have been designed for different purposes\footnote{Once, one of my summer undergraduate research students calculated twin primes, and statistics about twin primes, using Excel (much to my surprise), which makes this point well.}). Thus if a computational system is ``universal'' then one should be able to use it to emulate any other computational system, and that includes humans and their deductive systems (what we defined above is meant to cover human activity too). Indeed Turing \cite{Tur1} writes,
 \begin{quote} 
 The idea behind digital computers may be explained by saying that these machines are intended to carry out any operations which could be done by a human...\footnote{He goes on to exhibit the sort of instructions you might give a computer via a ``domestic analogy'':
 \begin{quote}
 ``Mother wants Tommy to call at the cobbler's every morning on his way to school to see if her shoes are done, she can ask him afresh every morning. Alternatively she can stick up a notice once and for all in the hall which he will see when he leaves for school and which tells him to call for the shoes, and also to destroy the notice when he comes back if he has the shoes with him.''
  \end{quote} }
   \end{quote}
He then writes  
\begin{quote}
Digital computers ... have been constructed, according to the principles we have described, and ... they can ... mimic the actions of a human computer very closely.
 \end{quote}
 Surely Thurston, Voevodsk\u{\i} and others were aware that   this interpretation of  the Church-Turing thesis means that there is no hope that a different calculation method will lead to a  \emph{theoretically} better way, and therefore  all we have left is what we do  \emph{in practice}.  It is important to note that Turing's remarks (discussed at the end of the last paragraph) do not contradict ``universality'' since computers can each answer everything that is answerable, but not questions that are not within its design parameters (i.e. its axiomatic framework).

\subsection*{A suitable timeframe} 

Can there be a   mathematical statement that is provable within our standard axiomatic framework, for which every proof is   too long for humans' timespan?\footnote{In practice, computing resources are finite. If every atom of the universe were a computer, working at light speed, from now until the end of time, that would still be a finite number of steps, since the number of atoms, the speed of light, and the remaining time are all finite.}
In practice, mathematical problems often come in families denoted $A_1,A_2,\dots$ (like ``factor the integer $n$''). Such a family
is    in the complexity class P if there is  a polynomial time algorithm\footnote{That is, the algorithm takes 
no more than $n^c$ steps, for some $c>0$, for every $n\geq 2$.} to resolve each $A_n$, and is in the complexity class NP if there is a polynomial length proof for each $A_n$. What is the difference? For NP we don't require an algorithm to find the proof, only that it exists and so P$\subset$NP.\footnote{We did not require that we know how to find the shortest proof, just that there is one that is so short. Thus ``NP'' stands for \textbf{N}on-deterministic \textbf{P}olynomial time, where   ``\emph{non-deterministic}'' means that one can not necessarily find  that proof in a short time. We believe that P$\ne$NP, the outstanding open question in theoretical computer science.}  For example, no one knows a fast algorithm to factor integers, but given the factors it is easy to multiply them together to verify a factorization, and so factoring is in NP.

Simple counting arguments show that the set of NP-families of provable correct statements is a rather tiny subset of the set of all families of provable correct statements. Thus in practice, no machine, no algorithm, no human or computer can hope to access all, or even a tiny fraction of, the  correct and provable mathematical statements deducible from a plausible axiomatic system.

In \cite{LoWig}, Wigderson states that
``problems in NP are really all the problems we ... mathematicians, can ever hope to solve, because  [we need to] know if we have solved them\footnote{The same can be said, for example, for physicists positing the existence of fundamental particles, in that they need to propose an experiment that could falsify or verify their existence.}... this is the very definition of NP: a problem is in NP exactly if you can check if the solution you got is correct.''

``If P$=$NP [then] all these problems have an efficient algorithm, so they can be solved very quickly on a computer.... if P$=$NP  then everything we are trying to do can be done. [which] is why P$=$NP   ... would be so consequential. However, I think most people believe that 
P$\ne$NP.''
 
\section{Formal proof vs  culturally appropriate, intuitive explanation}

Venkatesh \cite{Venk} states that ``a proof is generally understood to mean an argument compelling consensus'' indeed that 
 \begin{quote} 
 A proof is  defined by the fact that [it] should
  induce uniform agreement about [its] validity, without any need for replication.
 \end{quote}
De Toffoli's definition \cite{DeT2} is that ``a mathematical proof is a correct deductive argument for a mathematical conclusion from accepted premises that is shareable and verifiable a priori'', moreover that 
 \begin{quote}
 It is the public availability of arguments that allows other mathematicians to perform quality control. This is essential to filtering out non-proofs ... Shar[ing] mathematical arguments ... [is] a necessary condition of mathematical justification.
 \end{quote}
 However if proofs are  judged by the community, then can a proof be said to be objectively valid?
 Indeed De Toffoli \cite{DeT1} believes that 
 ``criteria of acceptability for rigorous proofs are not carved in stone ... but  are indexed to a mathematical community in a particular time''.
 
Even more, human proof verifiers are most satisfied when, after examining a proof, they can reproduce it, or something like it, in their own words.  This way the new knowledge is not just part of the sum of all human knowledge but rather it  is contextualized, and  known of its own accord as part of much larger picture.
  In Plato's \emph{Theaetetus} (148e-151d) Socrates makes the important point  that   a failure of rigour can lead   to a 
  \begin{quote}  miscarriage of ... thoughts [so]  ... the lack of understanding   becomes apparent to themselves and to everyone else. \end{quote}
  I would claim this supports the positive consequences of deducing new ideas from old and contextualizing them, rather than considering the new ideas in abstract isolation.

But these considerations do not address ``objective proof''.
Surely we want objective proof? Ideas that have been verified so that they work in any context.
 If so,  how  can we achieve that?
By making the proof and its understanding part of a larger program of understanding, or by ensuring the proof follows from the axioms, what we might call a ``formal proof''.
Nathanson \cite{Nath08} writes for many in claiming:  
  \begin{quote}
 How do we know that a proof is correct? By checking it, line by line.
  \end{quote}
Hales   \cite{Hal3} explains what formal proofs are and why some people find them appealing:
 \begin{quote} 
A formal proof is a proof in which every logical inference has been checked all the way back to the fundamental axioms of mathematics. All the intermediate logical steps are supplied, without exception. No appeal is made to intuition, even if the translation from intuition to logic is routine. 
 \end{quote}
 This goes a little beyond Hacking's \emph{Leibnizian ideal} \cite{Hac} of a proof:
  \begin{quote} 
Every step is meticulously laid out, and can be checked, line by line, in a mechanical way,
 \end{quote}
 in that a formal proof requires one   (at least in principle) to create a logical path deducing the claim from the appropriate axioms.
 
Perhaps the most common fear is that in reading a formal proof, line by line,  one can be  convinced by each step yet not perceive the whole, not see global error
 
 One might also ask how such a formal (or Leibnizian$^+$) proof would be verified? Do we just believe it on the say-so of the formal prover (human or otherwise)?
 If not, who would be the ``independent authority'' that does the verifying?  And what skills does that independent authority bring to the job of verification? Would they check it line-by-line, and perhaps miss the wood for the trees? Or would the authority have a strong intuitive notion of what is going on (and if so where do they get that intuition from)? Hales goes on to write
\begin{quote} 
Thus, a formal proof is less intuitive, and yet less susceptible to logical errors.
 \end{quote}
How does one interpret this claim? What errors are not ``logical errors''?  Why does the formal approach necessarily reduce the possibility of overlooking issues that arise, especially if the reader is not expected to see the big picture?

\subsection*{More logical fallacies about formal proof}

Formal proof tries to avoid ``intuition'' as being imprecise. There is a belief that  a presentation with only formal steps can help an independent authority more easily verify the proof. In 1895, Peano wrote
  \begin{quote}
 Imprecise ideas cannot be represented by symbols,
 \end{quote}
 which implies   that ideas  represented by symbols must be precise.
 Even if we agree with this extraordinary (and unproven) claim, it is easy to misinterpret it by believing that once  you have translated your mathematical problem into symbols you are guaranteed to not be wrong (rather than not be imprecise); of course, one can be precise and wrong!  Let's suppose that we are careful about meaning. Then who is the ``independent authority'' that does the verifying?  And if that authority has limited understanding   and little intuition, could a subtle error slip by?\footnote{Many errors can occur, not only mathematical, but also, for example, in the interpretation of symbols and language.  It is easy to invent an unambiguous protocol for any given identified issue, but can we know we have thought out all possible situations that need a protocol to avoid ambiguity or misunderstanding?}  
 So Peano suggested a formalization that necessitates identifying the fundamental mathematical ideas in an explanation and then finding a way to express these within a limited language (\cite[section 2.3]{Slim}). Moreover Peano wished to reduce proofs in his formal language (via various ``identities'') to as few symbols as possible to conclude that the mathematics has then been better understood.  Even if we agree that the mathematics is better understood in this context, what happens when one  translates Peano's shortest proofs into the standard lexicon so more humans can understand? Will this be the most desirable proof? Will it help us to see the next results?

Moreover, can a formal proof verifier  see how to reproduce the proof on their own?  Is it really verified, if it is not easily reproduced? It is then part of the union of human knowledge but is it really known of its own accord?

Formal proofs typically chase the details of  a proof back to the axioms. It is like a child tirelessly asking ``Why?''  (until one gets back to immutable truths), but at the end of that process, does the child remember what they asked at the start and how they got to the end?   A  proof like this is little better than the answer ``Because I said so'', no matter who is the  objective trusted authority. When we are  functioning participants in a community we expect   answers we can understand, interpret, appreciate, and even use if possible. We can also be excited to find an alternative or  clearer proof, though that plays no role in a formal system.


\subsection*{What other kinds of proof?}
 Eugenia Cheng \cite{Chen}  writes that philosophers\footnote{I think she means to write ``some philosophers'' and indeed some mathematicians.} believe that 
 \begin{quote} 
Thanks to the notion of `proof', we have an utterly rigorous way of knowing what is and isn't true in mathematics, 
 \end{quote}
 but  that  mathematicians  perceive formal proof as an over-focus on  precision.
So if we are so skeptical about formal proof, what works instead?
There is a wide gulf that separates traditional proof from formal proof.
Hales \cite{Hal3} writes, 
  \begin{quote}
Traditional mathematical proofs are written in a way to make them easily understood by mathematicians. Routine logical steps are omitted. An enormous amount of context is assumed on the part of the reader ... A trained mathematician [can translate] those intuitive arguments into a more rigorous argument.
 \end{quote}
 That is, a substantial amount of \emph{tacit knowledge} comes into understanding and contextualizing traditional proofs. This is reflected Hacking's \emph{Cartesian ideal} \cite{Hac} of a proof:
 \begin{quote}
After some reflection and study, one totally understands [the proof] and can get [it] in one's mind.\footnote{Hacking seems to suggest that to be satisfactory all details of a proof should be understandable at the same time. Although this is certainly desirable, I am not sure it is required for believing in a given proof, so long as one has at some time felt confident of each part of the proof, and of how the different parts join together.}
 \end{quote}
 But this does not quite explain what leads us to believe that a proof is correct.
 Cheng \cite{Chen} remarks, 
\begin{quote}
Although proof is what supposedly establishes the undeniable truth of a piece of mathematics, proof doesn't actually convince mathematicians of that truth ... Something else does.\footnote{Indeed  in a letter  to Dedekind in June 1877, Cantor exhibits a bijective map from Euclidean 2-space to Euclidean 3-space, commenting ``I see it, but I don't believe it''.}
 \end{quote}
  Research mathematicians therefore write ``proofs'' that are a convincing argument which perhaps could be turned into a formal proof, written primarily to enhance knowledge and understanding, while maintaining some level of rigour that convinces the reader that more could be done in that direction.  Mancosu  \cite{Man1} notes that this type of proof
\begin{quote}
does not bear directly upon some of the traditional foundational concerns, such as certainty, which have dominated much of philosophy of mathematics. 
 \end{quote}
 Nonetheless, we believe that such a proof, especially if it is widely understood, is ``robust'' and so less prone to error.    This has long been the ``community standard'': 
 \begin{quote}
 There is no ... mathematician so expert in his science as to place entire confidence in his proof immediately on his discovery of it...Every time he runs over his proofs, his confidence increases; but still more by the approbation of his friends; and is rais'd to its utmost perfection by the universal assent and applauses of the learned world.  \\ . \hfill ---  David Hume (1739)
 \end{quote}
 Hume is claiming that confidence in a proof stems from its robust nature under enquiry from technically competent, interested peers, not from some abstract verification. This system is not perfect. Significant errors are sometimes found in published proofs that have been  accepted to be  true,\footnote{And we would guess that there are many more errors out there, yet to be found.} and yet we stick with our system; it must have its advantages.
 
 \subsection*{The usefulness of a good proof}
When reading a  proof, a research active mathematician wishes to add to their own intuition and scope, not simply  agree that the proof's argument is correct.
The reader is not passive. She  wants to understand, to  synthesize and to use the ideas in her own research:
 \begin{quote}
 Mathematics \dots [does not] reward passive consumption. Understanding a mathematical paper is like visualizing a building based on the architect's drawings: the text and formulas are only a blueprint that the reader must use to reconstruct the author's imaginary world in her mind. If she does that, however, then the best mathematical theories have the same breathtaking quality as the image of Paris folded on itself.\footnote{A reference to the dream-architect scene in the movie \emph{Inception}.} The experience can be both exhilarating and addictive.
  \\ . \hfill ---  Izabella Laba \cite{Laba}
 \end{quote}

Different people get different things out of a reading and therefore a new research article can inspire new ideas in hitherto unforeseen directions. Even the same person can, at different times,  get different things from  reading an article, our understandings do change over time, sometimes even how we approach the whole area. 

Rather than chase proofs back to the axioms, most readers  rely  on the published literature (a \emph{library} of reliable knowledge),  on what is   already known.
This means the reader asks  ``Why?'' a reasonable number of times to be satisfied of the correctness of a proof (unlike what might happen in a formal proof), at least if the reader has enough current knowledge, in that they have read and mostly understood all, or most of, the references on which the article is based.  
 
 To summarize,   proofs are   accepted by community standards. This means that they might be wrong since we don't expect the details to be incredibly carefully checked; so  what is the purpose of proofs?  Perhaps the plan of the proof is the primary thing, explaining what ideas are strung together to prove the proposed theorem, so the reader can learn, retain and re-use the ideas.\footnote{Cheng \cite{Chen} compares this to legality (the proof) and morality (understanding that the proof is correct in principle). It would be nice if proof and understanding were synonymous but they are not, and we have to appreciate and  accept how they match and how they differ.}
 The details are usually of less interest, especially as an experienced reader can reconstruct them, though sometimes with an undesirable amount of effort.  
\subsection*{Writing a proof} 

 A proof is an explanation to a particular audience. 
 \begin{quote}
 Mathematical explanations are context-dependent since ``in different contexts different features might be salient'' 
\\ . \hfill -- Mancosu  \cite{Man2} quoting Lange \cite{Lang} 
 \end{quote}
 An explanation for an expert in the field (the first readers of a research paper) is different from that for a novice student in the area, which is again different from someone not in the area. Here we are talking about the level of detail required, and explanation of how arguments fit together. Today, the usual protocol is that papers are  written so that the abstract and half the introduction is accessible to a broad audience, an overview in the introduction to specialists (including novices), and the details for experts.
 
 Explanation can be given by examples and test cases, by analogy with already understood results, or by proving the result in some special cases that highlight the main ideas. A ``good proof'' to a mathematician is one that explains as well as proves, in fact the better the explanation the less danger there is in omitting cumbersome details. However there is a wider community that may wish to read and appreciate a proof of a result and then the details may become more useful.

Cheng \cite{Chen} claims that the purpose of mathematical communication is to turn the author's beliefs, via a proof, into a \emph{believed truth} of her reader. One cannot convincingly just state the belief; the colleague will only start to consider it to be true if there is plausible reasoning attached. But what is reasoning?  Will our correspondent believe the reasons?  Without the rigour and structure of   proof, without fitting the reasons into an appropriate framework, the reader will probably remain skeptical, or at least is less likely to remain skeptical with a well reasoned argument in the form of a proof. This is because a carefully worded proof helps allay the fear of ambiguity or misdirection.

\subsection*{Mathematical truth is beauty, and beauty truth} 
What makes a piece of mathematics ``feel right''? Cheng \cite{Chen} claims it is about what ``ought to be'', and not to be confused with ``useful, fun, intriguing, beautiful, proved in detail''.  This perspective helps understand what motivates the approach that many take to proving theorems.

The prolific and influential mathematician Paul Erd\"os claimed that an objective supreme being has a ``Book'' which contains the perfect proof for every true theorem, each of which is short and elegant. Short, so it is easy to verify, and elegant so one knows that the statement fits so well that it must be true. This is a wonderful conceit of professional pure mathematicians, when a concept is ripe to be understood then it should willingly yield its most succulent fruits. 

 The great Grothendieck did not   believe in big steps or   examples to elucidate progress: he believed that when a subject is ready the theory should be clear so there is no point in trying to push progress too fast.\footnote{To my taste, a sad reductionist theory of progress.}  Case-by-case analyses  may be complicated and clumsy but can help the researcher understand the patterns that  may lead
to a proof that makes a case-by-case analysis unnecessary. Proofs develop over time, and different authors choose to publish at different stages in the potential process.\footnote{This is one of many arguments against ``citation indices''. Half-baked good ideas will get improved by lots of people and so will be often quoted. A brilliantly thought out breakthrough, with a beautiful proof, may be difficult to improve, even marginally, and so get  referred to far less often.}
Some research mathematicians are loathe to publish anything but a ``final proof'' in which one can see why a concept is true in one  fell swoop.\footnote{Including Gauss whose  ``motto'' was ``\emph{Pauca sed matura}''
(``\emph{Few but ripe}'' in English).}
 Others are less selective, publishing a less complete theory, but partial results can encourage others and enhance the sense of community. It takes all sorts.


\subsection*{The language of research articles}
To a large extent mathematicians follow the library model of building upon past papers. A theorem is a key new result of the research paper.  Since research papers build on the library of knowledge, they typically quote what is known, and perhaps make simple modifications of what is already known to fit their needs -- these are the ``lemmas''. If such a simple deduction seems to be interesting in its own right,  then it might be called a ``proposition''. There are no hard and fast rules for the use of these words but one will find few variants.


One can argue that no two people view language in the same way, no matter how well defined. Indeed we get used to words in common usage, like ``epistemic'' seemingly used for a rainbow of connected but subtly different meanings.  But language is how we communicate proof, and it is rarely a perfect tool for that. Instead we hope that our description gives enough details that the reader is able to ``figure out'' or to  re-construct what we meant from what we have written, perhaps in their own way which may be subtly different (and that perhaps leads to something new). This argues for formal proof and language but we have seen that that has its pitfalls.

\section{What is an accepted proof in pure mathematics?}

Hales \cite{Hal3} writes that 
``philosophers tell us that mathematics consists of analytic truths, free of all imperfection'' and proceeds to ``prove that 
$1+1=2$.'' He then writes
\begin{quote} 
If only all proofs were so simple. Mathematical error is as old as mathematics itself. Euclid's very first proposition asks, \emph{on a given straight line to construct an equilateral triangle}. Euclid's construction makes the implicit assumption (which is not justified by the axioms) that two circles, each passing through the other's center, must intersect. We revere Euclid, not because he got everything right, but because he set us on the right path.
\end{quote}

The central idea of accepted proof is simple enough: Starting from agreed upon axioms we construct a proof of   given statements. Famously, Russell and Whitehead showed that if our axioms are consistent and logically independent then many seemingly simple statements 
take an inordinate amount of proving. So to advance far in mathematics we need to avoid going back to the axioms all the time. We need to build a    library of statements that we know to be true and are unambiguously stated.
Traditionally this library is stored in research articles, and synthesized in books. 

Then a researcher can advance the sum of all human knowledge by making logical deductions from what they can quote from our library.  By now there are millions of articles, and so there is a lot of scope for errors to creep in to the system. Articles might be flawed, or quoted statements might be misinterpreted (or even mistyped).  

\subsection*{The refereeing process}
To try to ensure that articles are correct,  we don't take the researcher's word for it that they have made correct logical deductions.  Good journals assure quality control by assigning one or several anonymous ``expert reviewers'' to a submitted article. These reviewers/referees usually begin by   judging whether what is claimed is true, new and interesting enough to be worthy of publication in that journal. Famously, Littlewood would ask\footnote{Ben Green, the managing editor at the \emph{Mathematical Proceedings of the Cambridge Philosophical Society} also asks his editors to \emph{persuade him} that an article that they recommend is indeed ``interesting''.}
\begin{quote}
Is it new?   Is it correct?   Is it surprising?\footnote{Venkatesh  \cite{Venk} argues that ``the value we assign to a work of mathematics is purely subjective, in the sense that it depends solely on the perception of that work, and not on any objective quality''.}
\end{quote}
 Then the referee is supposed to pick the text apart looking for logical fallacies. It might be that what is claimed is true but there are mistakes, possibly rectifiable, or that   explanations given are hard to interpret or somewhat ambiguous. The reviewer might not only identify errors and ask for further explanation, but might suggest significant changes (that they think will be improvements) and so get involved in helping the article develop to its final form. The idea is that, in this way,  the ``literature'' is interesting and safe to quote and believe for future authors.

\subsection*{The dangers of expert reviewers}
Referees are usually chosen because they are expert in the topic of the submitted article; often they have been referenced by the author.  Some referees are defensive of their work and are very picky when reviewing work by others on their favourite topic. Other reviewers are happy to see other researchers participate to explore the same questions. The outlook of those reviewers obviously affects how they review and what they find acceptable.  As an editor of several journals I look for reviewers who are keen, generous champions of their field but hold themselves to high standards, so when they review someone else's work they expect things to be done as well as possible but are willing to help some authors to perhaps more fully think through their arguments so that what is finally published does constitute progress in the field.

Any sentient reader will see the pitfalls in the process as I have just described it. Prejudices can enter the picture in any number of ways and editors must trust the reviewers to be fair, though occasionally something is obviously amiss in a 
review.\footnote{Though rarely.}  There are  relatively few people in esoteric research fields,\footnote{Venkatesh \cite{Venk} writes, ``The size and complexity of modern mathematics means that most papers are almost incomprehensible to us; our opinion of them can then only repeat that of others. The only people who can be involved in the formation of opinion about a given paper or a given question are those who interact with it in some way.  Now, the set of people who study the details of any argument themselves is very small; a much larger group acquire, instead, an awareness of its relationship to other existing work.''}
 and so all the experts probably know each other, or at least of each other, and regularly meet at conferences. They have many of the same research goals and so are often either collaborators or competitors (or both). There is little hope that they can be truly neutral when assessing each others' work.

So what are the alternatives?  Perhaps non-experts can  impartially judge submitted articles? After all, any   submission should build logically on the literature, and if written clearly the reader should be able to simply follow the argument. This assumes a lot. For example how easy is it to read and appreciate what is there in the published literature?  The latest breakthrough might quote several major advances from the last few years, so   to  follow the arguments of the new submission the referee must be able to understand and appreciate what was cutting edge not so long ago. In fact, to referee a top quality work, the referee will have to attain a fluency in the top works of the field in the recent past, and so become, in part, an expert herself.\footnote{And therefore refereeing can only be done properly by experts -- this feels like a  G\"odel-inspired argument!} For example, I was recently asked to referee a 104 page paper by a top journal, the authors building on their own  and others' work, amounting to several hundred further pages all published in top journals in the last   five years. Although I am an expert in the field it took me weeks to do an adequate job, even though I have enough confidence in the authors' technical skills that I do not feel that I need to check every detail.\  
  As an expert this is an appropriate use of my time (so as to keep up with the latest in my field of research), but one could hardly expect someone who has not immersed themselves in these questions and techniques to give up weeks of their time and perhaps find themselves incapable of doing an adequate job.

  \subsection*{Checking details}
The reader might flinch when reading that there are situations in which ``I do not feel that I need to check every detail'' when refereeing. From a pure perspective that is self-evidently problematic -- how can one review an article for errors and feel it appropriate to not check details? But this is where expertise comes in. If one knows the area then there are often proofs that are more-or-less ``standard'', so the expert reviewer can see that it looks roughly correct, and perhaps can identify quickly what part of the proof might differ from previous similar proofs in the literature, and focus their efforts there. Everything is new to the non-expert reviewer so they must check every detail, as they have little idea what precisely to scrutinise for where potential errors might arise.  Moreover, the expert reviewer is better equipped to help the author fix a faulty argument since they might have been stuck on the same issue in the past.

Mistakes commonly arise from  mis-applying the literature, for example by quoting a result out of context:
Most articles are written to be read from beginning to end, but   in the pages before stating a  result an  author might have included an assumption to be used throughout the article.\footnote{For example, in my area, analytic number theory, one might write ``Throughout we let $f(x)$ be a Schwarz-class function'' rather than repeating that in the  statement of every Lemma and Theorem.}
Then the researcher quotes that result verbatim, neglecting the earlier hypothesis and so misapplies the result. The non-expert reviewer might verify the statement of the quoted result but is unlikely to read the whole referenced article, and thus is unlikely to be able  to identify this issue, whereas an expert might have more feel for when a stated result is liable to be valid.
\bigskip

I hope we now agree that expertise has its advantages in refereeing, and is indeed unavoidable in practice,  although refereeing remains onerous no matter what. 
I have largely been discussing ``cutting edge work'' above, the sorts of advances that  truly make a difference. However the majority of submissions are not game-changers, but rather, small advances, and perhaps more easily reviewed. Nonetheless the same caveats hold, just to a lesser degree.  But here the expert can truly make the difference in helping an author who has good ideas but perhaps has not yet developed the technical skills to take their ideas all the way. In my experience as an editor many referees are encouraging and helpful in these circumstances, particularly if the author is an ``unknown'', explaining how they might modify what they have done to make the argument correct, or the theorem stronger or more general.\footnote{This is why I am against ``double blind-reviewing'' in which the authors' name is concealed from the referee, as it discourages the referee's generosity: It seems referees tend to assume an anonymous author ``should know better'' than to make that mistake, and so give a terse explanation about a mistake, rather than a helpful one.}

Finally if only experts are capable of reviewing big advances in cutting edge research what happens when a new submission contains a genuinely new idea? Something really different from what has gone before?  If all the potential referees are reading from the same hymn book, aren't they all likely to be skeptical of the validity of the new tune? Particularly if it is not-so-well explained or contains   non-lethal technical flaws? Certainly there may be initial (understandable) skepticism but the current system often works well. For example, when 57-year old unknown mathematician Yitang Zhang claimed in 2013  to have made an extraordinary breakthrough on gaps between primes using the deepest of ideas in new and surprising ways, the established community rapidly acknowledged that this was a research development of the highest calibre.\footnote{Initially the experts were hesitant about immersing themselves in such a difficult manuscript by an unknown, but once the experts started reading they quickly realized that here was something special.} It also helps that different journals act independently: for example, the paper in which Mordell  proved Poincar\'e's conjecture on the finite generation of rational points on elliptic curves, and posed his famous conjecture that would eventually be resolved by Faltings, was initially rejected by the London Mathematical Society \cite{Mord2}, presumably because it was not then a fashionable subject, but subsequently accepted by the Cambridge Philosophical Society \cite{Mord1}.

\subsection*{When building a theory that is not interesting}
In \cite{RTV} the authors write ``To judge the originality of ... work on the basis of a conception of the `existing body of knowledge' which comprises both secret and possible knowledge is intellectually callous''.  Certainly if one presents it like that.

Some authors come into conflict with the arbitration process by failing to take account of the culture of the community. For example, in combinatorics one learns proofs and then often adapts them to new circumstances (which may require substantial new inputs); few general theorems are of value since they will often cover only the known cases and uninteresting generalizations, and not be usable for new and genuinely interesting cases.\footnote{Gowers \cite{Gow} explains the value in such a perspective: ``While the structure is less obvious than it is in many other subjects, it is there in the form of somewhat vague general statements that allow proofs to be condensed in the mind, and therefore more easily memorized and more easily transmitted to others.''}    In his last few years Serge Lang  was upset that journal after journal independently rejected his general work of this nature in analytic number theory.
His perspective was  that if his precise result  had not been explicitly stated in the literature then no referee had the right to say that it was ``known'' (an axiomatic viewpoint), whereas practitioners felt that the proposed generalizations added little of value (a community viewpoint). There is  an old joke in refereeing:
\begin{quote}
What is new is not interesting, and what is interesting is not new.
\end{quote}
 In my experience, if this applies to a submission then the author may have difficulty understanding the report, no matter how well justified.

\bigskip

Whose responsibility is the correctness of an argument in an accepted (for publication) paper? The author, the reviewer, the editor, the journal?
It is generally agreed to always be the author's responsibility, no matter how much scrutiny the paper has gone through.
It is not the referee's job to assure the paper is correct but rather that they are supposed to have made a good, serious effort to verify the details.

\subsection*{The robust nature of proof}

We would like a system in which mathematical statements that are known as ``proven'', are reliably correct and always true, because a proof exhibits a theorem as the logical consequence of steps that are each verified. The proof should surely guarantee the theorem, so a competent mathematician does not need to verify each proof she uses herself (though that would be preferable). However this is not really feasible with the current system of proof verification of expert peer-reviewing.

We have discussed the advantages of expert peer-reviewing. It largely rests on the idea that familiarity simplifies scrutiny (at some risk of favouring conformity) so that papers are more quickly and accurately refereed. We mentioned the idea that a referee, having seen a certain type of argument before, does not feel it necessary to review every detail, perhaps only  focusing on those details that might be most likely, in the experience of the referee, to cause concern.  This supports a belief in the \emph{robust} nature of proof, at least at a high level of research. We believe that  not much   can go wrong with well-used technical tools and so we are prepared to make assumptions about what needs verifying,  usually true but perhaps not always. And then even if there is a mistake, experience shows that a simple modification should be enough to make the argument work. Any experienced researcher   does this regularly in developing their own work, so when they encounter minor technical flaws in the work of others, they tend to believe that they are fixable. Indeed many referees write in their report, a sentence like: 
``I do not need to see a revision as I am sure that the author will be able to deal with these issues''.

\section{Mistakes}

Our expert peer-reviewing proof-verification system can and does go wrong. Big results are widely used and  
errors   are sometimes found after publication. Usually not from someone checking the details even more carefully, but rather by another researcher  applying the claimed result and deducing something that is wrong, or implausible or at least unexpected and therefore suspicious. 
Let's discuss some notable examples.

\subsection{Goldston and Yildirim} In 2003 Dan Goldston and Cem Yildirim announced an extraordinary advance on what could then be proven about the shortest gaps between primes: We believe the \emph{twin prime conjecture} that there are infinitely many pairs of primes that differ by 2, but at that time it was not even known that there were infinitely many gaps less than a tenth of the average gap, and   Goldston and   Yildirim claimed they could get that down to any positive multiple of the average, no matter how small.  They distributed a preprint in which the novel part of the argument was written out in easily verifiable detail, though the proofs of some of the technical lemmas were sketched. Nonetheless those technical lemmas all had their roots in similar statements in the literature and so were believed; it was felt that Goldston and   Yildirim would just fill in the details later, the proof seemed robust and their work was widely believed, and prominent mathematicians told the scientific press what a great breakthrough it was, etc.  

Asking them for more details on two of the lemmas, they responded that one of the lemmas was proving stubborn but that they expected to sort it out soon and then would post a completely proved revision. That seemed plausible, so  Soundararajan and I then proceeded to further develop the 
Goldston-Yildirim method and, to our surprise, managed to deduce that there are infinitely many gaps between primes of size $\leq 16$, almost the twin prime conjecture! However something did not feel right and when we carefully traced our steps we found that   that stubborn lemma seemed to imply too much, which led us  to  a counterexample.\footnote{For the reader who knows complex analysis, the key issue revolved around how to move contours when integrating in very high dimension. There are aspects of high dimensional geometry that are different from the dimensions we are familiar with. The Goldston-Yildirim lemma worked   in low dimension but went awry in high enough dimension.}
This story again highlights the issues in accepting  proofs based on robustness but there is no obvious way to proceed otherwise.

This story though has a happy ending. Working with Janos Pintz, Goldston and Yildirim were able to find a correct version of what they had done, which eventually led to Yitang Zhang's breakthrough, and the Maynard-Tao modifications of the Goldston-Pintz-Yildirim method, and ultimately to the knowledge that there are infinitely many pairs of primes that differ by no more than 246, one of the great results of mathematical history.\footnote{Moreover there is a second happy ending: Green and Tao tailored the original Goldston-Yildirim sieve weights to complete their proof that there are arbitrarily long arithmetic progressions of primes.}

 \subsection{Wiles}
 In 1993 Andrew Wiles announced in a series of three lectures at a small conference at the Isaac Newton Institute in Cambridge, that he had proved Fermat's Last Theorem (FLT), using a very deep strategy suggested by work of Taniyama, Shimura, Weil, Frey, Ribet and Serre, to show that certain families of elliptic curves can be parametrized as had been conjectured (the \emph{modularity conjecture}).  Wiles had asked the conference organizer, John Coates, whether he could give three one hour lectures, refusing to tell Coates why but suggesting that it would be worthwhile. I was at the conference and rumours flew high. What was clear at the time was that Wiles presented a well-tuned sketch of his proof, even though having just three hours to explain the conclusions after seven years of thoughts was not enough. Wiles had the audacity to create a brilliant program of research to finish FLT from where Frey, Ribet and Serre had left it. Every step was either plausible with the then current technology, or he had a plan to make it plausible, and at least sketched the details. We were all  aware that filling in those steps might require some changes in the precise details of the plan, but he was confident that what he had was sufficiently robust that it could withstand a few missteps that might need correction. My impression was that he gave the talks a little before he was 100\% ready, but this was a conference that would be attended by most of the experts in the world, in a brand new research institute at Wiles' Alma Mater, and he felt the participants would prefer to hear about it now, since they could confirm that the steps were all doable or at least plausible (and even perhaps make useful suggestions).\footnote{In my memory one participant remarked that Wiles used a result that was known to be false but that something similar should be true and provable.   (I tried to verify this with that person  but he did not recall.)}
 
   As is well-known, when Wiles went on to fill in all those details, one part of the proof started to feel more shaky and stubborn no matter what his efforts, and he felt that he needed to retract his claim of FLT, at least temporarily.    Why had the community been so prepared to accept Wiles' claims, before his work had been released and refereed?  After all, anyone can claim such a result.\footnote{Indeed many have. Before Wiles' Theorem many of us used to be inundated with amateur claims to have proven FLT. Afterwards we were even more inundated, at least for a while, but felt less guilty about refusing to read them over.}
   The answer is about community culture. Wiles had already several major results that had changed the way we thought about algebraic number theory; he was already a professor at Princeton, and his every new paper already produced excitement. Moreover he was not known to have made any serious mistakes in his earlier work.\footnote{Actually I don't know if he had made any mistakes, but that is not relevant. I included the word ``serious'' to emphasize that no one worried too much about minor mistakes as long as the game plan of the proof was coherent and feasible, particularly for the experts.}
His description of his attack on (part of) the modularity conjecture was feasible to the experts in the room, indeed they could not wait to get their hands on a manuscript.\footnote{Ken Ribet, gave the lecture right after Wiles' third in which Wiles had announced FLT.  Ribet's lecture began with an excited room of mathematicians not quite knowing what to do with themselves. Ribet started by laughing a little nervously and said, ``After that, even I don't care what I am about to say''.}

When Wiles was forced to retract his claim, he worked to find a way to ``patch up'' his argument. First with comparatively minor re-thinking of that part of the argument, and subsequently with more dramatic changes, but he continued to work with the same overall plan. He knew his plan to be robust, indeed several parts of his proof plan have developed into interesting mini-subjects in their own right, that something should work to fill in the gaps, it was just a question of finding a workable path.\footnote{That path required further and deeper understanding than the original.   In retrospect, he had made it more difficult for himself by announcing his proof publicly, as he then felt under   pressure to get it done, whereas he had previously spent seven years working well away from the spotlight.}
While Wiles worked on the ``fix'' (eventually collaborating with his former student, Richard Taylor) there were different reactions from leading researchers in the  community. A few competitive souls told the press that what Wiles had was far from a proof, but the majority  seemed to be willing Wiles on; praising him for what he had already accomplished, acknowledging the brilliance of the approach he had created and expressing hope publicly that he would complete his Herculean task.\footnote{There were also some who felt that Wiles should publish what he had and let others ``have a go'' at fixing the hole.}

\subsection{Biss}
In the early 2000s, Daniel Biss made some important advances in homotopy theory, publishing papers in top journals like \emph{Annals of Mathematics} (indeed his work  developed ideas of leading mathematician and Annals editor, Robert MacPherson).  Biss's proof underwent the tough scrutiny you would expect from a journal like Annals. And yet there were mistakes, in fact the key mistake was submerged in the paper, in Proposition 4.5. Biss's result had seemed  ``feasible'' to the experts yet the mistake invalidated his main result. Similar mistakes appeared in other Biss papers in    top journals like \emph{Inventiones} and  \emph{Advances in Mathematics}, very likely refereed by  other top experts.  So what went wrong with the system here?  Do all the experts think so similarly, and accept what seems plausible so readily that such a mistake slips by, even in the most reputable of journals?  It gets worse though. Other Biss papers have mistakes; for example, in 2017 \emph{Topology and its Applications} ``retracted'' Biss's 2002 paper as ``the definitions in the paper are ambiguous and most results are false''.

 It is not wholly uncommon in some of the more abstract areas of pure mathematics that such issues arise with definitions or the simplest arguments (since they are often the least scrutinized). The flow of the overall argument is beautiful and persuasive so even if  some of the details have not all been ironed out by the time the author has gone public, that is usually a technical issue that will be rapidly resolved. Biss's  mistakes were identified by Nikolai Mnev in St Petersburg who quietly notified the author and other experts and  expected   Biss to rapidly resolve the issues or publish a retraction.  Biss worked with colleagues who believed that these issues might be resolved and yet were frustrated. Time passed and Mnev, after four years, felt he should notify a broader community (by publishing his observations on the arxiv), as he was worried other people would develop Biss's work without realizing that it might be wrong.  Four years can seem a long time to publish a correction or an erratum, but it might take a while in very deep and difficult research to identify the correct path forward, particularly under pressure.   
 
 In these last two examples the community held off for a long time  passing judgement on a mistake until things are clarified. There is perhaps no correct answer as to when we should   abandon a brilliant and seemingly robust proof plan.
 
 \section{Rethinking axioms and  language}
 
 \subsection{Mochizuki's rethink}
 Shinichi Mochizuki  is a highly respected researcher in arithmetic geometry who has made fundamental contributions to anabelian geometry. In 2012 he announced a proof of the abc-conjecture, perhaps the most fundamental open question in Diophantine arithmetic.  The proof stems from the creation of 
 ``inter-universal Teichm\" uller theory'' (IUTT), which he claims is a massive re-think of algebraic geometric aspects of Diophantine analysis.  Mochizuki did something very different by discarding a lot of the development of arithmetic geometry and replacing it by his new theory that he believes is more ``fit for purpose''.  Mochizuki did not compromise in his creation, the language is new, the concepts are, he claims, new and different, and he suggests that they are not translatable into the usual language and structures used in  arithmetic geometry.  Mochizuki has every right to  create his own theory involving a new language and structures. As long as his proofs can be shown to follow from appropriate axioms in the usual manner via some path, then they are correct.\footnote{Indeed many top mathematicians (one associates Grothendieck in particular) have found that the correct language and definitions can change the complexion of deep problems from being far from what is doable to being something that follows inexorably from the  theory. Indeed Scholze \cite{Harris} recently wrote: ``We perceive mathematical nature through the lenses given by definitions, and it is critical that the definitions put the essential points into focus.''}

 However, a practical problem arose in refereeing his works; Mochzuki had gone so far afield that there were no other experts to decide upon the correctness, and those  who would usually seem best qualified to judge research in this general area were not prepared to do so  without collegial help.
 Mochizuki uses a new language, and his four preprints giving his claimed proof of the abc-conjecture are about 500 pages long.
  He makes no attempt to compare his ideas with the standard lexicon, or to give   motivating analogies with the usual arithmetic geometric objects. 
  Moreover he turned down requests for more standard explanations from those experts stating that since everything is clearly defined in his work, those other experts have a responsibility to study Mochizuki's work from first principles.  This would have meant putting their own research careers on hold for a significant time, though arguably they would have gained much in doing so if they eventually agreed there is high value in Mochizuki's proposed revolution.
  
  Mochizuki's proof establishes a form of the abc-conjecture that is a bit different from the usual conjecture. Vesselin Dimitrov then deduced a more applicable version of Mochizuki's result, and found counterexamples, so that the result originally claimed   by Mochizuki is provably incorrect.   Mochizuki, claiming his proof was robust even if there were some minor glitches,  then weakened one claim so that the eventual result could not be disproved by Dimitrov's example. 
 Although it is not unusual that long and complicated proofs need some minor revisions, in this case the extremely minor change led people to wonder whether there might be counterexamples to the new version that one does not yet know about --  what then? How many ``tweaks'' are acceptable?
  
  A small coterie of established researchers  studied Mochizuki's works and believed the ideas and logic to be plausible. In the process of these studies Mochizuki   continued to revise what he had done to make that process easier. However that still left the majority of top arithmetic geometers unwilling to immerse themselves in IUTT, and so the community could neither satisfactorily accept nor satisfactorily discard Mochizuki's claims.\footnote{And this is not the place to comment on its correctness, since our focus is on how the community has functioned in this peculiar and complicated situation.}
  
  Peter Scholze, one of the top mathematicians in the world,  took the responsibility of  trying to translate   enough of Mochizuki's manuscripts into more common mathematical language so that he could judge for himself.  Believing they understood Mochizuki's strategy, Scholze and Jakob Stix identified a point (Corollary 3.12) at which they could not see how to proceed, and went to Kyoto in 2018 to discuss this with Mochizuki. Mochizuki met with them but claimed that their translations involved invalid oversimplifications of 
 IUTT so that their objections were invalid. However he was unable to persuade them that they were missing the point. They published a report \cite{SS} in which Scholze and Stix claimed that ``small modifications will not rescue the proof strategy''.
 
 This is  a situation in which the usual system has broken down because the  majority of mathematicians expect to be presented work that is presented as an addition to  what came before and following certain traditions. Not only in the writing but also in the communication. Mathematicians expect to be able to ask each other questions about  their latest work, either at conferences or by email, and expect clarifications in the most convenient language for all concerned.\footnote{For example,  Nick Katz at Princeton  writes with a great familiarity about certain algebraic-geometric objects in a way that I find difficult to appreciate. But in discussion he is prepared to discard all that and  to give as down-to-earth an explanation as possible, recognizing my limitations. This is arguably the ``community standard''.}

 \subsection{Voevodsk\u{\i}}  Vladimir Voevodsk\u{\i} was a great intuitive mathematician. Trusting in his intuition he  was   shaken to discover some of his work with Kapranov was wrong. Somehow his  intuition was insufficient and there were fatal mistakes in some of their arguments, and even counterexamples.\footnote{It took him 15 years from being made aware of the mistake by Carlos Simpson (via Kapranov), to acknowledging it since he wanted to find out where his intuition had gone wrong, and not worry about a detail that might have been correctable.}  Voevodsk\u{\i} highlighted in \cite{Rob15} what can go wrong with our expert-peer review system: ``A technical argument by a trusted author [like him], which is hard to check and looks similar to arguments known to be correct, is hardly ever checked in detail.''
 
 His mistake led him to propose a re-think of the foundations of his field: The \emph{univalence foundation program} launched by Voevodsk\u{\i}   attempts to reconcile Martin-L\" of-style dependent type theory with the traditional mathematical treatment of proofs and categorical constructions. This involves a lot of rethinking of foundations with a new approach and new language. However in this case, one of the main architects of this program, Dan Grayson, has gone to great pains \cite{Gray} to try to make the fundamental ideas accessible, and even translatable for other mathematicians, as one can surmise from the   article's title, ``\emph{An introduction to univalent foundations for mathematicians}''.  We will discuss this new program of study again a little later.
 
  Voevodsk\u{\i} also suggested a plan for the proof of the general Bloch--Kato conjecture in 1996. This provided a road map to the eventual proof in 2011. Along the way parts of Voevodsk\u{\i}'s plan were proved and parts discarded and replaced (most notably by Chuck Weibel and Markus Rost). However the entire strategy and plan were robust enough to find a way around the obstacles (albeit involving the creativity of a number of remarkable mathematicians). This is a big advantage of a community approach, led by a plan.

\section{Computers and proofs}

 Frenkel is quoted  in   \cite{Rob15} as saying, 
 \begin{quote}
 As the towering mathematical enterprise reaches new heights, with more intricately and esoterically turreted peaks, each practitioner drilling down into minutiae that could never be understood by colleagues even in the same field, computers will be necessary for advancing the enterprise, and authenticating the integrity of the structure --- a structure so complex in detail and so massive in scale that adjudication by a mere human referee will no longer be feasible.
  \end{quote}

 There are three main uses of computers in proofs:\smallskip
 
  $\bullet$ They can be used for calculations in establishing a proof;
  
  $\bullet$ They can be used to assist in verifying the logic of an author's arguments, perhaps interactively. These are called ``computer-assisted proofs'';
  
   $\bullet$ They can be used to prove claimed theorems, so-called  ``computer-generated proofs''.
   \smallskip
   
   We now briefly discuss these in reverse order.
   
 \subsection*{Computer generated proofs} They are in their infancy.  Common views include: 
    \begin{quote}
    I don't believe in a proof done by a computer ... I believe in a proof if I understand it. \hfill -- Pierre Deligne 
  \end{quote}
   \begin{quote}
I'm not interested in a proof by computer ... I prefer to think\newline
\.  \hfill  -- John H. Conway \cite{Rob15}
     \end{quote}
     \begin{quote} 
     Computer-generated proofs...teach us very little\footnote{And goes on to say ``If we go toward computer-generated proofs then we lose all the good that there is in mathematics --- mathematics as a spiritual discipline, mathematics as something which helps to form a pure mind.''}   -- Vladimir Voevodsk\u{\i}   \cite{Rob15}.
  \end{quote}
Since these top researchers do not want to believe a claim to be true
just because the computer says so, we might counter these remarks by insisting that a computer-generated proof should be designed to be human-readable, and to attempt to \emph{explain} the ``why'' as well as the ``how''. We will return to this viewpoint later in discussing recent developments.

Kahle \cite{Kah2} claims that once a computer finds and describes a proof it is relatively easy to verify it is correct, even if that is tedious (so a good job for a computer). However finding a proof is hard. The idea of ``trying every possibility'' suffers from well-known complexity problems but there are strategies to improve upon this, and we will discuss a different approach of 
Ganesalingam and  Gowers.

\subsection*{Computer assisted proofs}
Russell and Whitehead exhibited, in their arduous task, that for anything beyond a very trivial result, the number of logical inferences in a formal proof is too large to be adaptable, that the whole quest seems to only be of interest, in of itself. So it was not too sad when G\" odel put paid to their original purpose. However formal proofs are making a comeback!  No longer are they ``roped-off museum pieces to be silently appreciated, but not handled directly'' (\cite{Hal3}). Now we have computers that have the memory space to handle the length of more-or-less any proof, and (we hope) the logical resources to ensure that no steps are omitted (though this raises   questions about trust in computer calculations, as discussed below).

In the future proof-verification might   employ ``computer-assisted proofs'', since the author could interactively explain her proof to an appropriately designed proof-checker.  Indeed in 2008 Harrison \cite{Harr} at Intel, wrote that one of his goals for formalization is
\begin{quote}
Supplementing, or even partly replacing, the process of peer review for mainstream mathematical papers with an objective and mechanizable criterion for the correctness of proofs.
\end{quote}
Even now this would involve an inordinate drain on the author's time, but it might nonetheless be  useful in situations where there are many new definitions that must be correct, for example in the works of Biss and Voevodsk\u{\i} that we discussed above.

We  must surely be wary of believing in computer verified proofs for the old reason that we are translating mathematics into a specialized language.\footnote{In fact proof verification software is a spin-off of hardware verification software, and software now can formally verify 
that   high level computer languages  or microprocessors (or anything in-between) operate as claimed. For example Leroy \cite{Ler} created a formally verified compiler for the C programming language.}
To justify this one can  resuscitate Peano's  belief that in an appropriately designed language one can eliminate mistranslations and obtain ``precision'' so as to eliminate mistakes. However,  even if this is possible and you have achieved this dream, how would you prove that you have succeeded?

In Voevodsk\u{\i}'s work on his univalence foundation program he felt he no longer trusted himself on the details and designed a proof assistant. In \cite{Rob15} he claims that a proof-assistant can keep you honest: you lay down the plan, it builds the boring details. When it can't you have to refine your plan further;   like working with a mythical pedantic enthusiastic colleague. Roberts \cite{Rob15} writes that Voevodsk\u{\i} was  
\begin{quote}
jousting with the computer. He instructs it to try this, try that
\end{quote}
 rather like playing a video game.
 
    \subsection*{Calculations as part of proofs}
When computers were first able to produce inordinate amounts of data and information that might be used in a proof, there was some struggle to decide what to require of an author for the reader to  believe in the proof.  The consensus now is that the author must produce a coherent explanation and justification of the calculations, and describe them so that the interested reader can easily repeat them independently.\footnote{As self-evident as this might seem now, one can find ``philosophical discussions'' of computer verification in the literature, both by research mathematicians and philosophers, which expresses discomfort in accepting such calculations, even with appropriate protocols.  It is always initially difficult to accept technological change: For example, in 1994, a few years after the world wide web became ubiquitous, there was an MSRI workshop on  
``the Future of Mathematical Communication'' to discuss its potential. Like now, the internet contained a lot of unedited nonsense and many of the more senior  participants at the workshop seemed skeptical that one could create an online journal which maintained the standards of print journals.  Even when the obvious solution was proposed (using referees), several remained unconvinced.}  
  
  \section{Uses of computers in major theorems} 
  
  Most proofs include calculations by the author. Indeed authors routinely claim properties of the objects they are working with, without feeling a need to justify them in detail, since anyone in the field will know or at least expect these things. These are not definitions but claims that the object worked with is an example of something familiar and therefore obeys the properties we expect of such objects.
  There is no precise rule about  the amount of detail that an author needs to give: It depends on what she expects of her readers, and what she  believes her readers expect of her!  For example, one might assert that a given polynomial has no rational roots without further explanation, or one might state ``by Gauss's rational root criterion'', or one might work through the finite number of possibilities that Gauss's criterion yields. Sometimes the amount of calculation required to justify such a claim is substantial and so the author may feel compelled to explain how the calculations are done, so as to be helpful for a reader wishing to check the results, rather than give the data, which might not be enlightening at all.\footnote{Especially when there may be a large number of cases to verify for a simple property.}

  We begin by discussing examples where extraordinary calculations have been part of proofs, and then how proof assistants have subsequently helped.
  
 \subsection*{Short gaps between primes} In Maynard's work on small gaps between primes \cite{May} he needs to construct a sieve with certain properties which he showed follows from constructing a special polynomial of degree $\leq 11$ in $105$ variables, so that some function of this polynomial is $>4$.    To verify Maynard's proof one can follow his steps precisely, or one can veer off his path at some stage and see if one can get to the same conclusion via an alternative route.\footnote{Having spent many years thinking about similar questions, I was skeptical when Maynard first told me he had got the constant $>4$. I not only got the same answer proceeding as Maynard suggested, but also tried modifications since I worried that his persuasive writing had concealed a mistake.} By explaining things well, he thus gave the reader the opportunity to verify his game plan, and to infer that his proof is robust. The details of the calculation would have been less helpful than a description of the thinking behind each step.

\subsection*{God's number is 20} This refers to the number of twists needed to resolve Rubik's cube from an arbitrary starting position. That is, no matter how Rubik's cube is jumbled when you start you can get it back to the starting position in 20 moves, and there are positions that require 20 moves (as shown by Michael Reid in 1995). This  claim was made by Morley Davidson et al  in \cite{DDRK} in 2010. There are about $4.3\times 10^{19}$ possible positions of Rubik's cube.  With a little group theory the collaborators reduced this to resolving about 56 million positions in 20 moves or less, and then used 35 CPU-years of idle computer time (donated by Google) to nail these down.  Finding the 20 moves is a much bigger proposition than following the 20 moves already written down, so this is not so hard to check, and check independently, using their output file.\footnote{In practice a typical PC might be able to do $10^{15}$ basic computer steps (like adding two digits) in a day, which means that checking $<10^8$ positions should be easy in a day, but more than $10^{14}$ impractical in a year.}

\subsection*{The Classification of Finite Groups} 
This is another example of an overarching, robust plan, followed by an enormous number of parts to be filled in.  In this case even the plan was and is of such extraordinary depth and complexity that only a few seemed to be confident of it (though this has improved over time).  Filling in the steps was always going to be a massive job, with many participants and all sorts of room for error.   However here many of the steps are interesting in their own right so that more people were willing to be part of this great project. Advances meant that new structures emerged and the plan got modified to explain more while becoming technically simpler. One senior participant, Michael Aschbacher \cite{Asch1, Asch2} honestly remarks that at times he has believed that the classification has been complete, at other times not, and he is always certain that there are minor errors throughout. He believed in the robustness of the plan, that a full classification would be found along the lines claimed, albeit with some modifications, sometimes substantial. As yet the classification has not been computer verified.

\subsection*{The Four Colour Theorem (4CT)} 
 There is a lot of (philosophy) literature expressing epistemic dismay at the original proof of  4CT, some of which misunderstood what happens in the purported proof(s), so let  me try to clarify. The key  notion is, again, ``robustness''.
In all proofs  the  idea is to  show  that if a counterexample exists to 4CT then one can ``reduce'' it to find a ``minimal counterexample'' belonging to some finite computable set, and then to prove that no such minimal counterexamples exist since they would have to have too special properties (or be further reducible).\footnote{This is not an uncommon strategy in  graph theory.} This has, to date, been  complicated. 

In Appel and Haken's 1976 original proof of 4CT, their discussion of how to reduce is not entirely well organized and in parts is difficult to verify. They needed to show that 1478 different subgraphs can be reduced further; they used a computer but without convincing documentation. Moreover the computer algorithm, as described, is complicated (involving 487 steps) and difficult to verify as valid and non-self-contradictory. Nonetheless Appel and Haken did create a convincing plan for proving 4CT in which one could see the point of each major step and why it should work. Here the devil was in the details, and even if these were implemented correctly they were always going to be hard to verify, believe and build upon.\footnote{I have heard it said that quite a few mistakes were identified but each could be dealt with by relatively minor modifications of the details; however there were so many such fixes that few people in the area had faith in the details of the end-product.} There is a robustness to the overall plan but the implementation was unconvincing.

In 1996, four leading graph theorists, Robertson, Sanders, Seymour and Thomas \cite{RSST} (RSST) decide to rework the existing proof to make it  more believable. They followed the same overarching plan as Appel and Haken but looked for simplifications in the implementation. They had to computer reduce just 633 subgraphs (still a large number, but the authors went to great pains to make their construction transparent). More important is that their computer algorithm only involves 32 steps, so that a strong mathematician in the area could spend a day or two and believe that the authors had successfully covered all the options.  This new proof was consciously written to be easily verifiable (though still long and complicated) and no serious mathematician doubted that 4CT is proved, at least by community standards.  That is,  this RSST proof is  sufficiently robust that one expects to be able to easily patch up any misunderstanding that some future researcher might unearth.

A few years later 
Gonthier \cite{Gon}  verified the RSST proof using the  Coq v7.3.1 proof assistant, by developing a ``formal proof'' that covers both the 32 steps and the 633 subgraphs that needed to be reduced; the proof then ``depends [only] on the correct operation of several computer hardware and software components'' and is not specific to this proof, which feels more robust.
Most importantly many of the programs used here had been used in other calculations, which makes one feel that no bug that is specific to what is done here could have crept in. Moreover Coq can   produce a proof ``witness'' (albeit human-unreadable).\footnote{A \emph{witness} is a relatively short verification that a problem has been correctly solved. For example, to prove to you that I have factored $147573952589676412927$, I can simply produce the factors $193707721 \times 761838257287$; you do not need to repeat the steps that led me to these factors. To verify you can simply multiply the two factors together.}

\subsection*{3-dimensional sphere packing}

Hales' 2005 proof \cite{Hal} of the Kepler conjecture on packing unit balls in 3-dimensions\footnote{In 1606 Sir Walter Raleigh asked a mathematical acquaintance to help determine how many cannonballs he could stack  on board a ship; by 1611 Kepler conjectured, in his paper, ``On the six-cornered snowflake'', that the hexagonal pattern (and a closely related alternative) used by greengrocers to stack round fruit would be the   best way to squeeze as many spheres as possible into a large space.}
 was assisted by large-scale calculations.
 Like in the proofs of 4CT, Hales' strategy
was to reduce a minimal counterexample to Kepler's conjecture  to some finite computable set of possible arrangements of spheres 
which would have to satisfy some extraordinary constraints. 
A large calculation  could then enumerate these possibilities and rule them out.
However there are so many cases and the refereeing was extremely onerous leading Lagarias \cite{Lag} to write
\begin{quote}
The nature of this proof . . . makes it hard for humans to check every step reliably. . . . [D]etailed checking of many specific assertions found them to be essentially correct in every case ... [This] produced in these reviewers a strong degree of conviction of the essential correctness of this proof approach, and that the reduction method led to nonlinear programming problems of tractable size,
\end{quote}
a conditional but appropriate endorsement.
A second  proof \cite{Hal2} in 2017 used the   HOL Light and Isabelle proof assistants. The purpose of each step is fully described  in
\cite{Hal2} but the researcher certainly still needs to take a lot on trust.

\section{Computer error} 

One of my biggest concerns about computer proof systems is people's tendency to assume that once a program has been well-implemented it is reliable, and then to mistake ``reliable'' for ``free from error''.

We have all   dealt with computer systems which supposedly  ``never make mistakes'', yet they make them, whether it be the computer system for your credit card, phone company, bank,  airline or your home university. The reason could be programming or implementation errors,  a non-understanding of the possibility of your particular situation or even a hardware problem.
 Are these issues   avoidable? Can we correct these programs and computers to be trouble free?
 
 \subsection*{Dealing with all situations that can possibly occur} Turing \cite{Tur1} wrote:  
 \begin{quote}
 It is not possible to produce a set of rules ... to describe what [to] do in every conceivable set of circumstances.  One might for instance have a rule that one is to stop if  one sees a red traffic light, and to go if one sees a
green one, but what if by some fault both appear together? One may
perhaps decide that it is safest to stop. But some further difficulty may well arise from this decision later. To attempt to provide rules of conduct to cover every eventuality, even those arising from traffic lights, appears to be impossible.
\end{quote}
And then, even if a program is impressively accurate and reliable, what about   when upgrades appear? If we avoid upgrades,\footnote{Which has the disadvantage that one cannot take advantage of new developments.
In 1990 Donald Knuth made the decision to never again upgrade TeX, no new features, only bug fixes; most scientists prefer to use newer typesetting systems like LaTeX, based on Tex, which are regularly upgraded.}
 can we eventually perfect a system? And even if a system is ``perfect'' in that it really does respond appropriately to all situations that can arise, how would we know? (ie How do we know we have listed all feasible situations? How do you prove such a thing?) 

 \subsection*{Computer hardware reliability} 
The gold standard for computer hardware and software in the non-academic world is that there are no longer any complaints,\footnote{In other words, outputs that are obviously inconsistent with other information.} but this does not translate to guaranteed proof verification!
Former Intel President Andy Grove    said 
\begin{quote}
No microprocessor is ever perfect; they just come closer to perfection
\end{quote} 
in discussing hardware reliability.
It is feasible that  a small mathematical  error in a widely used computer chip could be exploited to defeat widely used cryptographic protocols which would put
all e-security  at risk. Indeed, in 1993 Pentium released a chip which they subsequently found had a hardware bug   affecting its floating point processor. Rather than recall the flawed chips they decided to keep quiet and correct the problem in updates. In June 1994, Thomas Nicely discovered the error while computing a number theory constant as precisely as possible; and subsequently it was found that the error could be detected when dividing certain seven digit integers by each other in several different softwares that used floating point arithmetic.  Even then Pentium resisted a recall until IBM refused to ship their product.\footnote{Pentium claimed the  recall   cost them \$ 475 million. One can understand why they hesitated!} Shall we elaborate on how this story affects our discussion?  

--- Perhaps there are errors in chips today that are more obscure and so less likely to be detected. Moreover if we do not independently analyse the results of large scale calculations (as Nicely did by comparing the answers he got with those that were already known) then how will we spot such subtle errors?   

--- Are seemingly small errors in a chip's output really worth the cost of fixing for their manufacturers?  Moreover manufacturers will rarely reveal concerns about their products (so as not to put off potential purchasers), so even if they knew about a fault how likely is it that that information would be widely shared?

\subsection*{Computer software reliability} 
There are also many problems with computer programs.
Hales \cite{Hal3} notes that commercial software contains about 1 bug per hundred lines of code, and perhaps 1 per 10,000 lines when programs are super-focused on being careful.\footnote{Like in the space program.}
He writes that  
\begin{quote}
corporations ... keep critical bugs off the books to limit legal liability ...
only those bugs should be corrected that affect profit.
\end{quote} 
 Moreover  correcting bugs can be problematic since that process often creates new bugs.\footnote{It is one thing to make a minor change in what you are working on when it is fresh in your mind and you can see the big picture as well as the small details, but it is difficult to regain that perspective when coming back to something complicated years later to make a hopefully minor fix that one does not wish to spend a lot of time on.} Harrison \cite{Harr} remarks that there may be more mistakes in a physical engineering project but since the questions there are typically continuous in nature, a small error makes little difference, whereas computer engineering acts on problems that are discrete. Although this gives one the opportunity to find exact solutions, it also means that one minor error, one mistyped digit, is more likely to lead to a major fault.
\bigskip

Perhaps the only (sensible) answer proposed to  these particular computer problems is to repeat the same calculations on different computers run by different chips, using different software. It would then be highly unlikely to get an error at the same point in the ``proof'' on several different systems. If results match from two such systems that are sufficiently independent then it seems extraordinarily unlikely that there is a problem. Therefore we feel that, in studying proof verification, we can safely ignore these computer implementation problems, though it is not clear what this means for guaranteeing (i.e.~proving) that a proof is correct.


\section{Protocols for automated theorem checkers/provers}

What features would be desirable in a proof checker?   The proof verifications that we have discussed above (4CT and Kepler's conjecture) worked interactively with a person to construct a proof from the ground up.  Those proofs are not human readable, but they can attempt to confirm further claims made by humans, for example minor variants to improve our trust in them. Indeed the more they show that they provide understanding, by helping us go further, the more faith we might have in a computer proof.

However, is a  ``proof'' that cannot be understood in detail,   really a proof? We have discussed how the purpose of a proof is not just to establish truth, but also to enhance understanding. If it cannot do that then  what use is it to the community?  

Why should we trust the output from a proof verifier or a prover if we can't read them?
Can  proof-assistants be self-correcting if they can only checked by their own internal logic? Indeed, it may well be that they continue to propagate a subtle error. 

Surely these programs need a community to verify their proofs? Perhaps their output may be independently verified by using different programs; in effect, we propose refereeing computer proof verifiers output within their own community!  This assumes that they work in a common language which adds extra burden to the different designs.  In this way humans might believe a  computer-verified proof, via an independent computer verification, and so the proof verifier becomes a trusted, objective, expert authority (that is, a referee).  The computer programs will use their  community to obtain a worthwhile seal of approval. In this way we can design the future based on what already works.\footnote{Rather like how e-journals used the refereeing process to establish, in their new context, the integrity and  standards of traditional print journals.}

\subsection*{Proof presentation} 

To believe in a proof that goes back to axioms, we need  to have a common language so that the proof can be independently verified. We have seen that it is not feasible for a human to do the verification, and that humans gain little from this process, so how can we make it more useful?  It seems evident that a proof verifier could also output a human readable proof. It could learn the types of high-level arguments that humans understand and appreciate, and then present its proof not only for human verification but also to help enhance human understanding.

The automated prover might select results from the existing library of verified results to build a short, person-readable argument to deduce the latest advance and so fit into the well-established protocols of how the community agrees on the correctness of a proof.  Anyone may use any result that has been previously established. Although each step in each proof is computer verified, back to axioms, one hopes that as more researchers contribute to the system, interaction will move towards something resembling the high-level practice of mathematicians. A system like this is  user-friendly and should become an integral part of the mathematician's arsenal.

Typically theorem-provers are  interactive, the user being able to give it hints.  
The user enters statements into the proof-verifier,    based on simpler objects that the machine already knows about. 
The proof assistant will determine whether the statement is `obviously' true or false based on its current knowledge. If not, the user enters more details. The proof assistant therefore forces the user to explain their arguments in a rigorous way, and  to fill in simpler steps than human mathematicians might feel they need.\footnote{Users report that they often learn a lot during the input process!}
For example suppose a proof needs ten lemmas. Some the theorem-prover will see and resolve quickly. Others it might be stuck on and the user gives it more details until the computer can see its way to a proof. In so-doing the program learns more, and maintains a library and is perhaps more efficient when it next encounters similar issues.   

 For now proof  assistants can't read a textbook, they need it all to be broken down for them by  humans.\footnote{And this human/machine interface can lead to problems. Indeed some definitions input into Lean by different users  have been inconsistent -- those we know about   have been corrected. But sometimes there are ambiguities in the literature. For example Kevin Buzzard pointed out to me that \emph{topos} can mean different things to an algebraic geometer and   a logician and sometimes the difference can be quite subtle. There is no clear way to deal with this dichotomy.} Proof assistants can't judge whether a mathematical statement is  interesting or important,  only whether it is consistent with what it has been shown.  It should eventually require less help, perhaps much less help. We have no idea when  (and whether) it will be able to generate its own proofs.

\subsection*{The uncertainty principle of objective proof verification} The history of mathematical practice suggests that

\centerline{\emph{The less one questions a proof,  the more susceptible it is to error.}}

\noindent This  important principle strongly suggests  one must find a wide variety of ways to explain and to verify any given proof, even a computer proof, and to look at it from as many different perspectives as possible. 

\section{The future of proof}

Computer generated proofs give too many details (often by a factor of more than 100) to be human readable, whereas human generated proofs give too few details to be computer verified!
Hence, we have to agree on what constitutes a proof, including how much detail is required.
If we accept that the computer needs to work with human-style proofs\footnote{After all, isn't that our purpose in creating machines!?} then it will need to be able to fill in  details  to justify the missing steps; the key difficulty is for the computer to independently determine what that requires. A famous example is the following proof of the irrationality of $\sqrt{2}$:
\begin{quote}
If $\sqrt{2}$  is rational then we can write $\sqrt{2}=a/b$ where $a$ and $b$ are coprime integers, so that 
\[
 a^2=2b^2.
\]
Therefore $a$ is even and we can write $a=2A$ so that 
\[
b^2=2A^2
\]
and $b$ is even. But then $2$ divides both $a$ and $b$ which are coprime, a contradiction.
\end{quote}
In 2006, Wiedijk \cite{Wied} noted that no  computer-verifier could take this text  as input and verify it as correct.

 We find that we need to return to the question of what a proof is to help our proof checker. Avigad \cite{Avig} explains that  a proof is a communication which provides sufficient information to establish that the purported theorem is true. Beyond correctness, it can be evaluated with respect to    background and interest.
All of the major 2006 computer-verifiers  could  have constructed a proof of the irrationality of $\sqrt{2}$ based on the above argument, but these would have been  of overwhelming length, and no one but themselves could check their own proofs. Moreover, since their languages are so different they could not even  check each others' proofs. Thus, as Kahle \cite{Kah} notes, the issue is not only to produce human readable output but  to work out how the proof-verifier can appreciate the structure of an abstract mathematical proof, and  represent that for others to read and analyze. It also needs to  \emph{explain} the mathematics,  not necessarily give the shortest  proof, but to be understandable, convincing and potentially easily reproducible by the reader.

\subsection*{Robustness and fragility}
We have discussed how traditional intuitive mathematical proofs are robust in that minor flaws can be fixed, and so we can have confidence in them, even without certainty. On the other hand, formal proofs are remarkably fragile in that if we find any errors, it puts much more into doubt -- once a formal proof is in any way wrong it calls into question all sorts of aspects of its formulation, particularly with our inability to read the details and therefore thoroughly review its claims.

It is not uncommon to have several seemingly different proofs of the same theorem, and it is worth asking when two are really the ``same proof in disguise''. Is it when the key ideas are the same? Or should we expect more to be in common? It is unclear (see \cite{DeT1}). Moreover sometimes the disguise is quite convincing and one cannot easily recognize the common threads.
This leads one to ask whether a  formalization of a given intuitive proof is going to be the same proof? When an intuitive proof is dissected into what is required for, say, Lean to work with, it will look very different and rest on a rather different looking library of knowledge.  And how different will the same proof look when modified for a different language?  

Buzzard gives a couple of great examples (in private correspondence):

--  Leaner A\footnote{A \emph{leaner} is someone who implements a proof in Lean.}  might prove Pythagoras's  theorem by assuming that the triangle lives in $\R^2$, changing coordinates to make the right angle at the origin and then proves that 
 \[
 \int_{(0,a)}^{(b,0)} 1 ds =\sqrt{a^2+b^2}.
 \]  
 Leaner B, formalizing Euclid's actual proof,   takes an  abstract Euclidean plane, never defines area, but rather defines
 what it means for two shapes to have the same area. These are not just different proofs, these are different statements
 of Pythagoras's theorem, which emerge from using two   different (but valid) models of the Euclidean plane.
 
--  In Lean real numbers are defined to be equivalence classes of Cauchy sequences; but another prover might define real numbers using
 Dedekind cuts. So a proof of, say,    the mean value theorem will look very different in the two provers, even if they must be equivalent since they must both  rely on the fact that a non-empty bounded set of reals has a least upper bound. 

\smallskip

Formalizers  hope   that  formalizations might help us see what is necessary to use in a proof (like which axioms).
The idea is to create a ``proof compacter'' that somehow recognizes how to shorten proofs, compacting them into as small a space as possible. The human only participates once the proof compacter cannot do more.  The hope is that this ``shortest possible proof'' will have its advantages though I do not think it is evident  why the shortest proof should use a minimal set of axioms or necessarily be advantageous.

\subsection*{Can computers generate their own proofs?}
  ``Machine learning'' typically develops its understanding in  simple ways as a result of clever algorithms.
Creating a large database and analyzing it with specially formulated tools can be startlingly effective (like Google Translate or ChatGPT) but this is not the same as developing intuition (or even simulating intuition effectively).  There is a lot of money and a lot of publicity surrounding the subject of ``machine learning'' and some other forms of ``artificial intelligence'' but rest assured that many hyped advances are either  exaggerated or easily explained in terms of well-designed algorithms and extraordinary computing power.\footnote{The impact of the Google search engine is more-or-less uncorrelated to how simple (though elegant) the mathematical ideas are behind it.} There are as yet no ``thinking machines''.

It is difficult to know how we can move forward in this direction, as relatively little is understood about creativity and intuition, and how we move from one understanding to a rather different one. To simulate this on a machine seems very far away.\footnote{Indeed the ideas behind the latest exciting developments in machine learning do not portend any real understanding; see Melanie Mitchell's wonderful book \cite{Mit} for a forensic discussion of what underlies some of these developments.}
  Ada Lovelace (1815-1852), who posited the concept of (what we would now call) a computer program from Babbage's early calculation machine, and even supposed that a computer could be taught to compose music, wrote
\begin{quote}
[Babbage's] Analytical Engine has no pretensions to originate anything. It can do whatever we know how to order it to perform.
\end{quote}
 Before computers, librarians were often credited with knowing a lot more than they really did (as the gatekeepers of so much knowledge). Computers are  much bigger repositories for knowledge, more accessible and less proscribed by others, and can achieve some surprising feats; it is not surprising they get credited with powers they do not yet possess.

 Can we design a computer verifier to learn and think like a human?
Ganesalingam and   Gowers \cite{GG} embarked on designing an automated theorem prover that proceeds rather differently than stripping everything back to the axioms. Rather, the idea is to proceed like a human, to model the way humans think, to produce proofs that read like a human proof. In their design they require input and output that is understandable and uses the standard lexicon (to be ``user-friendly''),\footnote{It is considerably harder for a machine to understand such input than to produce such output, a bit like my experience with speaking different languages.} informative solutions that are not just verifications, and that further capability can be easily added to the program to incorporate new concepts, problem-solving techniques etc. In short one should be able to interact with the program as one does with a (marvellously retentive) human.  There have been previous human-oriented programs but the machine-oriented programs have seen much more progress over the last twenty years (see section 1.3 of \cite{GG} for some history).

Machine oriented search programs tend to try many options, which has the disadvantage of combinatorial explosion but with some \emph{tactics} (like pruning search trees) this can sometimes be well-managed.  Humans have to avoid severe combinatorial explosion so bring in their tactical awareness earlier in the process. For example if a theorem to be proved has several hypotheses and conclusions the theorem prover doesn't know which is most important (and so to be focused on) and basically tries the different statements in a random order to find connections; a human might quickly examine the important conclusion and start puzzling as to how one might get there from what is known.  A computer might try very similar substitutions over and over in its search tree whereas a human might see from one example that a particular type of substitution (say where $y$ is linear in $x$) cannot work in general. So there are strategies that humans have that are atypical of large search strategies. As yet we do not know how to list all these differences.

In \cite{GG} Ganesalingam and   Gowers  remark \begin{quote}
For the majority of proofs that mathematicians find, there is some kind of `story' to tell of the ideas that give rise to the proof. Typically, such a story will be a high-level overview of the main difficulty and how it is overcome, where `overcome' means that the problem is reduced to one or more problems where that difficulty no longer occurs. Often this reduction is achieved by means of a well-chosen intermediate statement that turns out to follow from the initial assumptions and imply the conclusion. The intermediate statement itself is typically found not by means of a brute-force search but by a process of approximation: one might make a guess, find that it is unhelpful, understand why it is unhelpful, and use that understanding to guide the search for a better intermediate statement. These characteristically human techniques enable mathematicians to penetrate deep into `proof space', but the set of proofs that can be discovered in this way forms a tiny fraction of that space. It seems almost a truism that human methods will be useful for programs that want to find these special proofs that human mathematicians are so mysteriously good at finding.
\end{quote}

Humans are pretty good at selecting which technique or techniques to try and a new one can improve their efficiency at proving things. On the other hand a new technique simply expands a computer's search space and might well decrease its efficiency.  
 
 In terms of describing human proof adequately we are  still at the stage that we don't know how to describe what a good proof is, but we know it when we see it. The work of both humans and machines can be defined in terms of ``tactics''.
 For now we can best use  interactive systems based on ``tactics'' that are  designed to mimic human reasoning.
However it seems to be difficult to get machines to recognize which tactics to use when, that is to provide an order in which to try  different tactics and perhaps to adjust the future tactics or their order depending on intuition gained from  a tactic that has  just failed. In \cite{GG} they ask whether to let the computer learn from its past experience, to try to devise a theory that better mimics human choices, or work with a mix of the two.  It is also hard to decide when humans use certain tactics, like proof by contradiction.

 Eventually we will need to play Turing's ``imitation game'' \cite{Tur1} with machine created proofs; that is, their proofs should be indistinguishable from great human proofs. So we set the machine up against a Fields' medalist, ask them both a tough but doable question, and see whether we can determine whose proof is whose.\footnote{In the imitation game \cite{Tur1} an interrogator tries to distinguish between  a human who tries to prove she is a human and    a computer who tries to fool the interrogator into believing it is a woman. (In the original, pre-Turing, game the part of the computer is played by a mischievous male.) This is now known as the ``Turing test''.} Turing \cite{Tur1} notes that since  computers are universal (via the Church-Turing thesis) they can  perfectly imitate anything that can be computed including human interaction, so inevitably this will be doable (and indeed computers will eventually be able to perfectly imitate each other).\footnote{He also made some predictions, for example that by 2000 the imitation game will have been mostly won by computers, and that by then ``one will be able to speak of machines thinking without expecting to be contradicted''. He was over-optimistic but, on the other hand, he did get a lot of things correct!}  
 
 In  \cite{GG}, Ganesalingam and   Gowers selected problems to prove and got thousands of independent readers to try to distinguish which proofs were by their program and which by real people  (see https://gowers.wordpress.com/2013/04/14).
 The results are encouraging, though of course this is not the Turing test since it is not an independent arbiter that selected the problems.
 
 It is hard to predict the future here. Things are moving fast and brilliant people are getting involved.
Kevin Buzzard (see the next section) takes the view that 
\begin{quote}
The more people are familiar with the software, the sooner interesting things will happen
\end{quote}
which is a compelling perspective.  

\section{The Lean Theorem prover}

The Lean Theorem prover\footnote{http://leanprover.github.io/} has taken the research mathematical world ``by storm'' in the last year or two.
Lean ``is a functional programming language that makes it easy to write correct and maintainable code. \emph{You can also use Lean as an interactive theorem prover.}'' It has been made popular in pure mathematics by notable Imperial College arithmetic geometry  professor Kevin Buzzard
(who reports on his progress on the blog \cite{Buzz}).
He writes 
\begin{quote}
I believe that digitising mathematics is important, for the simple reason that digitising anything enables you to do new things with it.\footnote{And goes on to write ``Currently the computer proof systems we have are not good enough to tell mathematical researchers anything new about the $p$-adic Langlands program [on which Buzzard is an expert] or other trendy Fields Medally things, so the top mathematical researchers tend not to be interested in them. However [this does not mean] that they will never be useful to us, and the point of [this] project is to make it happen sooner.''}
\end{quote}
  Buzzard's initial focus was to follow the axiomatizing dream by attempting to show that all  the undergraduate syllabus in pure mathematics can be justified (partly interactively) in Lean (a bit like moving through   levels of a computer game).  Thus many of the technicalities could be verified by Lean, leaving the student to  higher level thinking.  Buzzard refrains from predicting too much of what is to come from this new capability.

Buzzard is well aware that this creates the possibility of an enormous shift in the proof culture of mathematics. He involved undergraduates in creating proofs in  Lean for aspects of the undergraduate mathematical canon, but confessed (to me) that he was unsure how much of the proofs   these undergraduates understood (though ``perhaps they understood other things that are at least as valuable''). An independent project\footnote{https://www.researchgate.net/project/Learning-about-proof-with-Lean} is investigating how differently these  undergraduates think about mathematical proof.

\subsection*{Research level mathematics}
Buzzard writes 
\begin{quote}
In the near future I believe that maybe computers will be able to help humans like myself (an arithmetic geometer) to do mathematics research, by filling in proofs of lemmas, and offering powerful search tools for theorems ... but there is still a huge amount of work to do before this happens.
\end{quote}
In early June 2021 this hope was realized when Lean verified the part of an argument that recent Fields' medallist, Peter Scholze, was unsure about
in his notes on Analytic Geometry \cite{Sch2} with Clausen. In a guest blog post  \cite{Sch} Scholze writes 
\begin{quote}
I find it absolutely insane that interactive proof assistants are now at the level that within a very reasonable time span they can formally verify difficult original research.
\end{quote}

Developing Lean's formal version of the proof involved interactions between Lean and Leaners (the people formalizing the proof in Lean), and between the  Leaners and Scholze. 
The Leaners input the Clausen-Scholze manuscript line-by-line  into Lean, which  created ``a clear formulation of the current goal'', and then the Leaners would refer back to the manuscript to figure out how to proceed with the next few steps. If necessary they would contact Scholze for clarifications.
 Scholze writes:\begin{quote}
  Sometimes it was then realized that even on paper it does not seem clear how to proceed, and the issue was brought to [my] attention ... where it was usually quickly resolved.
\end{quote}
Thus Scholze's intuition assisted Lean (and the Leaners), which never needed to look more than one or two steps ahead to follow and formalize
the Clausen-Scholze proof.
 
 Lean did pick up imprecisions, in particular that a certain infimum need not be a minimum (which had been assumed) and required some modification of the original proof. Scholze notes that 
 \begin{quote}
This was precisely the kind of oversight I was worried about when I asked for the formal verification.
\end{quote}
 More interestingly Scholze felt that he learned ``what actually made the proof work''. His interaction with Lean via the Leaners allowed him to make various parts of the proof  more explicit (and so elementary) and to evaluate certain key constants that he had wanted to better understand. 
 Gonthier \cite{Gon2}  found something similar back in 2005 while developing his formal proof (in Coq) of 4CT:
  \begin{quote}
Perhaps this is the most promising aspect of formal proof: it is not merely a method to make absolutely sure we have not made a mistake in a proof, but also a tool that shows us and compels us to understand why a proof works.
\end{quote}
Indeed in the process of his interaction with the Leaners (particularly Commelin)   Scholze realised that he could get away with a weaker theorem and therefore eliminate the  stable homotopy theory used  originally.

Scholze warns that he cannot read the formalized proof; allowing him to do so   is  an objective of the Lean developers.\footnote{Lean's formalization of the Clausen-Scholze proof contains tens of thousands of lines of code, about 20 times longer than the original.   It is not easily human readable.} 
There are obvious issues; for example it might be ``obvious'' for a human that $A=B$ from the definition of $A$ and $B$, yet it takes quite a while for Lean to concur.\footnote{Akin to the old joke in which a professor is asked whether some step in a proof really is obvious. The professor goes silent, desperately thinking, until after   half an hour he replies ``Yes, it is obvious''.} Nonetheless Scholze feels that this is a 
 \begin{quote}
 landmark achievement ... to take a research paper and ... explain lemma after lemma to a proof assistant, until [it has] formalized it all.
 \end{quote}

\subsection*{Are there epistemological advantages of Lean?}

Is there  reason to claim that Lean makes  less unjustified assumptions without realizing it (as long as they appear to be consistent) than humans?
It seems that the process is  similar to what happens when human mathematicians verify each other's work
(asking pedantic questions, wanting to know more about the definitions and to see some details filled in, which can lead to revisions). Scholze explained to me\footnote{In the comments on the blog.} that

\begin{quote}
It was exactly the interactions with the [Lean chat] that convinced me that a proper verification was going on. To me the [fact that it was a] computer didn't matter at all, it could also have been Ofer Gabber.\footnote{A mathematician who is known to insist on the right details.}  If I [had] simply got a blanket `This is all correct' stamp, I would have been extremely suspicious, as I was completely sure that I made some slips.
\end{quote}

Scholze had been nervous of a very complicated and technical Theorem, and had been aware of some nasty details that needed sorting out. He wrote 

\begin{quote}
The confidence that it is correct comes from ... seeing that during the process, [the Leaners] seemed to run into all the little nasty details that I expected (and sometimes [had] not [expected]).
\end{quote}  
He ended up by stressing that 
\begin{quote}
This whole experiment was a genuinely human experience, very similar to going through this with a very careful colleague.
\end{quote}

Scholze's description of the proofs themselves being (currently) presented obscurely by Lean is frustrating (so it can feel like one's less communicative colleagues), but there is no reason for that aspect of Lean to not be improved. 
For now we have ``no sense of the terrain'', only where the peaks are and whether we have reached them, so what is it that  can truly be said to have been learned?  Perhaps now that  Lean seems like a viable proof verification assistant, its design can find a more intuitive way to describe the proof it has constructed?  Often after a first human proof scopes out the terrain, the next proof gives a more enriching description of the ideas involved, so why not Lean? This would be useful and move mathematics forward whereas,
 \begin{quote}
 For now, I can't really see how [Lean] would help me in my creative work as a mathematician.\hfill Peter Scholze in \textsc{Nature}
 \end{quote}

 
\section{Myths of objectivity} 
In confirming that a proof is correct we believe that we can recognize and establish an objective truth. But can we? It is easy to 
believe in one's own objectivity, or that of an ``unbiased machine'', but are such beliefs valid, or are they self-serving? 
To help highlight our potential fallacy, we now recall two famous examples of so-called objectivity,  
(arguably) objective in their own time and their  original paradigm, but certainly not now.

\subsection{These self-evident truths}
For a long-time in our society biblical guidance was considered to be objective. Turing \cite{Tur1} remarks that 
in the time of Galileo, the quotations
\begin{quote}
The sun stood still . . . and delayed going down about a whole day\\
\mbox{} \hfill -- Joshua 10:13 \\
He laid the foundations of the earth, that it should not move at any time
\mbox{} \hfill -- Psalm 104:5
\end{quote}
were considered by many to be an objective refutation of the Copernican theory.

Social commentators like Donna Haraway  \cite{Hara}  explore the abuse of   belief in ``objectivity''.\footnote{For recent in-depth discussions on the relationship between objectivity in mathematics and social ontology, see the theme issue \cite{CT}.}
  For her primary issue, sexism, it is not hard to make the case: for so long  society has tolerated a belief that  one's gender or race quantifiably affects one's ability, rather than societal issues, whether they be social class, upbringing, opportunities, etc. Few today would argue that those earlier   beliefs in objectivity were anything other than self-serving 
 and most people today believe that one
should aim to replace those beliefs by a more objective understanding. However,   Haraway remarks that \emph{any} proposed objective standard is really a subjective ``power move, not a move towards ... truth''.  
Indeed to say something is ``objectively proven'' (that is, objective within some current paradigm) seems to quash any potential objections.

Poor scholarship continues to support this kind of prejudiced ``objectivity''. For example, in 2005
Larry Summers, while president of Harvard University, chose to try  to explain the low numbers of women in STEM jobs:  At that time some test scores of cognitive abilities for twelfth grade women showed less variation than for men, indicating  fewer women at the top end (as well as at the bottom end). Identifying factors that might cause these different statistics, he focussed on   ``issues of intrinsic aptitude'' (as well as ``\emph{lesser factors} involving socialization and continuing discrimination'').   His excuse for selecting   ``intrinsic aptitude'' to be  the most compelling was that he ``wanted to be provocative''.\footnote{Sadly, representative groups  censured him on the grounds of sexism, rather than   inadequate scholarship. This led to their protests being ignored when President Obama selected Summers to be President of the National Economic Council, thus continuing to add an air of respectability to his poor scholarship.} Today, and even then, one might view such provocation as ``problematic''.

To go beyond highlighting bias, and to not allow the issue to become bad people's ``bias'' versus good people's ``objectivity'', feminists 
\cite{Hara} claim that bias is ingrained into any social construct, arguing that there is a 
``collective historical subjectivity''.\footnote{Haraway also notes the desire to find a legitimate feminist ``objectivity'', but this is evidently paradoxical, in that one is in danger of repeating the same errors that one has been at pain to identify in others.}

Haraway \cite{Hara} highlights that ``parables about objectivity'' are told to mathematics  students during  their training, even though this is not what mathematicians do in practice.  She argues that there can be no truly ``trusted authorities'' but rather an earthwide network of connection, respectful of different perspectives. Indeed she claims that
\begin{quote}
 Science has been about a search for translation, convertibility, mobility of meanings, and universality,\footnote{All issues touched on in this article, and indeed part of what we have referred to as the ``community standard''.}
 \end{quote}
but one fails these criteria when work is only being compared to the hegemony.
In the context of this article, Haraway's work suggests that the community standards that have served us well in the past should illuminate the path forward in the rapidly developing age of computer proofs. That these standards, verifying that the plan and immediate details of a proof appear to be correct, are the closest we can get to objectivity within any given paradigm.

\subsection{Objectivity and infallibility}  Can there be infallible proofs? Proofs that are ``unconditionally and eternally'' correct?
Leibniz and Newton's infinitesimals were mostly accepted in the 17th century, yet were replaced by the $\epsilon-\delta$ proofs as derived by Cauchy through to Weierstrass. Will our current proofs of the basics of calculus withstand the test of eternity? Are they infallible? Will they adapt to all the questions that will be asked of them? How can we know? We have a framework within which they seem to be incontrovertible, but will that framework seem appropriate in the light of future understandings?

\begin{quote}
A triangle. This seems to be extremely simple, and you'd think we \dots know all about it \dots  Even if we prove that it possesses all the attributes we can conceive of, some other mathematician, perhaps 1000 years into the future, may detect further properties in it; so we'll never know for sure that we have grasped everything that there is to grasp about the triangle. And this holds also for bodies, for their extension, for everything!  \hfill -- Descartes (16/04/1648) in conversation with Burman.
\end{quote}

And what about objective truth within the current framework? From G\"odel we know we can't come close to verifying all truths   since  the set of true statements is far larger than the set of provable statements so we need to restrict our attention to what is provable.\footnote{Even more we need to restrict our attention to the (much smaller) family of statements in the complexity class NP, which loosely means those for which there is a proof or solution that can be verified on a human time scale,  a key issue in theoretical computer science. This was discussed with a few details at the end of section 2.}
Even so how do we know whether a given  proof is correct? Tarski \cite{Tarski} notes that ``intuitive evidence is far from being infallible, has no objective character, and often leads to serious errors'', and he views the  ``subsequent development of the axiomatic method ... as ... restrict[ing] the recourse to intuitive evidence''. This argument suggests that we can only hope to use a formal proof to obtain  sound verification of an objective truth, but 
can a proof verifier become infallible and so be the ultimate arbiter of what is correct?  Without doubt a well written program may uncover new problems in established proofs, and create proofs that are more difficult to challenge, but that is not quite the same as infallibility.

The mathematical philosopher Avigad \cite{Avig2} claims that 
\begin{quote}
According to the standard view, a mathematical statement is a theorem if and only if there is a formal derivation of that statement, or, more precisely, a suitable formal rendering thereof.\footnote{So, it seems that the Clausen-Scholze result only became a theorem that Avigad would recognize, when Lean asserted the proof to be so.} 
\end{quote}
An attractive definition of a theorem, but resting on the assumption of the infallibility of a ``formal derivation of that statement''.\footnote{Avigad goes on to claim, ``When a mathematical referee certifies a mathematical result, then, whether or not the referee recognizes it, the correctness of the judgement stands or falls with the existence of such a formal derivation,'' a claim to an extraordinary overview of mathematical process.}
Obtaining an intuitive proof and a formal proof (based on the intuitive proof) are attractive goals,  different ways of assuring that there are no easily recognizable errors, but it is a stretch to believe one can, in this way, assert objective truth.  Avigad's claims rest on a view of objective truth that does not reflect any consensus (despite Avigad's confident phrasing), but rather serves to justify the focus of a certain sub-community. Indeed this is endorsed by the extreme claim of his colleague Azzouni \cite{Azzou},
\begin{quote}
Formalized proofs have become the norms of mathematical practice.
\end{quote}
Not quite sure how this became Azzouni's ``norm'' as it is so far from true for the majority of the mathematical  research community. The day may come when formalized proofs play something other than a peripheral role, but that day has not yet arrived.

   Detlefsen  \cite{Detle} better understands mathematical practice:
\begin{quote}
Mathematical proofs are not ... generally presented in a way that makes their formalizations either apparent or routine. This notwithstanding, they are commonly presented in a way that does make their rigor clear ... at least by the time they're widely circulated among peers ... 
\end{quote}
That is, the author of a community-accepted rigorous proof rarely concerns herself with formalization, although the formalizer must surely be concerned about rigour (but it is then odd of Azzouni to suggest that the formalization makes the proof significantly more trustworthy).
 Avigad \cite{Avig2}  also argues for  formal proofs since  ``providing less information only exacerbates the problem: if even a complete presentation of a formal derivation cannot be checked reliably, providing strictly less information can hardly provide more confidence.''   So  how much explanation is enough?  Personally I prefer a clear one page proof, than a lengthy turgid treatise that fails to appreciate what is important even if it dots a few i's,``filling a  much-needed gap in the literature''!
 
 Returning to the theme of the subsection ``Robustness and fragility'' in section 11, it is interesting that formalizers feel their proofs are more trustworthy when they are evidently so much more fragile. Moreover if, as Avigad seems to claim, a  proof is not a proof until it
 can be  formalised then one must ask which formalization is the correct one? Different Leaners might produce very different formalizations of the proof of Pythagoras's Theorem so which should be the accepted one? Do proofs come in ``equivalence classes''
 (that is, proofs that are the same in disguise are in the same equivalence class)? If so, what are the criteria for deciding  which proof belongs to which class, and if there is more than one equivalence class then what does that infer about the fundamental nature of proof?  Should machines strive for a ``Book'' where only the best proofs are given?

\subsection{Deus ex machina\nopunct } literally means ``god from the machine''.  It describes  the viewpoint of some in the formalized proof community.
Believing in infallibility (of their own code) leads some programmers to not appreciate that their programming might occasionally be wrong!
That only by  communicating one's ideas will their ideas be accepted. Indeed in  \cite{BaWi} Barendregt and   Wiedijk confidently assert that a  putative proof is verified if ``the small number of logical rules are always observed'', and so once they believe they have done that they create an  unenticing 
 description of the formalization (see figure 8 of \cite{BaWi}). They need to appreciate that a proof is a social compact, and so make the extra effort to inform   interested readers.

\subsection{The Clay Millennium prizes} 
In 2000 the Clay Foundation announced a prize of a million dollars for the resolution of any of seven famous mathematical problems.
A solution can receive the prize only two years after it has been published in a refereed journal, and ``has achieved general acceptance in the global mathematics community''. These rules leave little doubt that the framers only have faith in community-based proof verification, and even then feel that it takes a while to be sure.

\section{Will machines change  accepted  proof?}
In this article I have asserted that proof verification does little to change the central tenets of proof as a social construction. Moreover that  there is little added value in learning that a program claims a proof has been verified (without providing more helpful information to increase the reader's understanding). Nonetheless we can expect that efforts will be made to make those formal proofs more accessible and hopefully useful (for example, they might eventually remove extraneous ideas from intuitive proofs). Indeed 
Patrick Massot \cite{Mass} recently announced that software tools are   being developed to automatically convert formal proofs into  human-readable interactive proofs, allowing  a reader to dig progressively deeper
until she reaches a claim that she believes with no further explanation, whether in a formal or intuitive proof.

Since the chess program Deep Blue defeated world chess champion Garry Kasparov in 1997, machine learning programs have become increasingly good at strategy in board games, by now easily beating the best in the world at both chess and Go.  Indeed in late 2017, DeepMind's AlphaZero was switched on, played only games against itself for 24 hours, and achieved a grandmaster level in both games (and in shogi), making it literally ``superhuman''.  It has a relatively shallow search tree (80,000 positions per second in chess compared to 70 million for some other top software), but compensates for that with better tactics.\footnote{Chess grandmasters found its play ``alien'' with ``insane attacking chess'', for example sacrificing a queen and bishop to exploit a positional advantage, something a human would be unlikely to do.} Can such  ``reinforcement learning''  lead to a greater depth and variety of mathematical proofs?
Might it create proofs that are more surprising at times than human ones, like its chess play?

In October 2022, DeepMind introduced\footnote{https://www.deepmind.com/blog/discovering-novel-algorithms-with-alphatensor}   AlphaTensor, which was set  the job of multiplying together two matrices of given dimensions as efficiently as possible; that is, with the least number of multiplications. AlphaTensor improved what had been known in about 20 cases \cite{Fetal} but never by more than 5\%, and  it is particularly striking that it does not appear to have deduced any new general theorems, nor indicated how new types of general theorems might be found.\footnote{In pure mathematics, finding extensive examples if often a precursor to better understanding, and thus truly new Theorems.} This indicates that this new technology can surpass human observation that comes directly through calculation in areas where there are few strong theoretical ideas, but there is no indication yet that DeepMind's algorithms will lead to the creation of new, deep theorems. Perhaps one day  it will be able to better identify the reasons that its tactics are so successful, which might help better understand what is going on.

No one has yet built a  quantum computer that can calculate faster than a classical computer even in specially selected questions.
Nonetheless the theory suggests that certain parallelizable problems may be \emph{much faster} to resolve on a quantum computer, most famously Shor's quantum factorization algorithm, which has caused vast resources to be pumped into ``post-quantum'' cryptography.\footnote{The security of several important cryptographic protocols are based on the difficulty of factoring so if factoring becomes easy then these protocols will no longer be secure. This worries powerful people so there is a lot of money going into research in this area.}  For us the question is whether quantum computing could be   adapted to the task of finding proofs (for example, if one uses a ridiculously large search tree).  Carlos Simpson pointed out to me that we might well run into the problem of obtaining a proof but having no idea how it was found, and maybe so long as to be uncheckable on any classical computer.

Raw computing power extends our mathematical capacities in many significant ways and, paraphrasing Kevin Buzzard, the more that people invest in the possibilities,  the sooner interesting things will happen.  Venkatesh \cite{Venk} suggests that mathematics might 
``be greatly altered; its central questions and values ... very different from those to which we are accustomed.''  This
``will enhance our ability to do mathematics but also will alter our understanding of what mathematics is''.
Indeed it is only a matter of time before we learn how to uncover tremendous possibilities for mathematics and for proofs revealed by computing power,  software, and brilliant programming ideas.

 \bibliographystyle{plain}

\begin{thebibliography}{99}



  
  
\bibitem{Asch1}   M. Aschbacher,
\emph{The status of the classification of the finite simple groups},
 Notices Amer Math Soc 51 (2004), 736--740

\bibitem{Asch2}   M. Aschbacher,
\emph{Highly complex proofs and implications of such proofs},
Philos Trans R Soc A 363 (2004), 2401--2406

\bibitem{Aus} J. Auslander,
\emph{On the roles of proof in mathematics}, in B. Gold \& R. A. Simons (Eds.), \emph{Proofs and other dilemmas: Mathematics and philosophy} (pp. 61-77). Washington, DC (2008) Mathematical Association of America.

\bibitem{Avig} 
J. Avigad,
\emph{Mathematical method and proof},
Synthese 153 (2006), 105--149.

 \bibitem{Avig2} 
J. Avigad,
\emph{Reliability of mathematical inference},
Synthese 198 (2020), 1--23.

\bibitem{Avig3} 
J. Avigad,
\emph{Varieties of mathematical understanding},
Bulletin AMS 59 (2022), 99--117.

\bibitem{Azzou}
J. Azzouni,
\emph{Why do informal proofs conform to formal norms?},
 Foundations of Science, 14 (2009), 9--26.
 
\bibitem{BaWi}  Henk Barendregt and Freek Wiedijk,
\emph{The Challenge of Computer Mathematics},
Phil. Trans. R. Soc. A. 363 (2005), 2351--2375.

\bibitem{Buzz} 
Kevin Buzzard
\emph{The Xena Project}, a blog at
https://xenaproject.wordpress.com/

\bibitem{CT}   Paola Cant\`u and Italo Testa (eds.)
\emph{Mathematical Practice and Social Ontology},
Topoi \textbf{42}  (2023), 187--344.

\bibitem{Chen}
Eugenia Cheng,
\emph{Mathematics, morally}, 2004 (preprint).


 
\bibitem{DDRK}
M. Davidson, J, Dethridge, H. Kociemba and T. Rokicki,
\emph{God's Number is 20}, https://www.cube20.org/.

\bibitem{Detle}
M. Detlefsen,
\emph{Proof: Its nature and significance} in ``Proof and other dilemmas: Mathematics and philosophy'' (B. Gold \& R. A. Simons, Eds.)
Washington, D.C.: Mathematical Association of America (2008), 3--32.

 

\bibitem{DeT1} Silvia De Toffoli, 
\emph{Reconciling Rigor and Intuition},
Erkenntnis (2020).

\bibitem{DeT2} Silvia De Toffoli, 
\emph{Groundwork for a fallibilist account of mathematics},
The Philosophical Quarterly 71  (2021), 823--844.


  \bibitem{LoWig} Bjorn Ian Dundas and Christian F. Skau,
  \emph{Abel Interview 2021: L\'aszl\'o Lov\'asz and Avi Wigderson},
  Notices of the AMS \textbf{69} (2022), 828--843


 
\bibitem{Eck}
M. Eckert, 
\emph{Water-art problems at Sanssouci -- Euler's involvement in practical hydrodynamics on the eve of ideal flow theory},
Physica D: Nonlinear Phenomena 237 (2008), 14--17.

\bibitem{Fetal}
Alhussein Fawzi, Matej Balog, Aja Huang, Thomas Hubert,
Bernardino Romera-Paredes, Mohammadamin Barekatain, Alexander Novikov,
Francisco J. R. Ruiz, Julian Schrittwieser, Grzegorz Swirszcz, David Silver, Demis Hassabis and Pushmeet Kohli, 
\emph{Discovering faster matrix multiplication algorithms with reinforcement learning},
 Nature 610 (2022),  47-63.


\bibitem{GG}
M. Ganesalingam and W.T. Gowers, 
\emph{A fully automatic problem solver with human-style output}
 J. Automat. Reason. 58 (2017),  253--291.
 
 \bibitem{God}
Kurt G\"odel,
\emph{The modern development of the foundations of mathematics in the light of philosophy},
draft of a 1961 presentation to the American Philosophical Society,
 Collected Works, Vol III, Oxford (1995),   375--388.

 \bibitem{GY}
 Dan Goldston and Cem Yildirim, 
\emph{Small gaps between consecutive primes},
(retracted preprint)

\bibitem{Gon}
Georges Gonthier, 
 \emph{A computer-checked proof of the Four Colour Theorem} (2005), preprint.

 \bibitem{Gon2}   Georges Gonthier,
\emph{Formal proof -- The four-color theorem},
 Notices Amer Math Soc 55 (2008), 1382--1393.


 \bibitem{Gow}   W.T.~Gowers,
\emph{The Two Cultures of Mathematics},
 in ``Mathematics: Frontiers and Perspectives?,
Amer Math Soc, 2000.





\bibitem{Gray}
Daniel R. Grayson, 
\emph{An introduction to univalent foundations for mathematicians},
Bull. Amer. Math. Soc.   55 (2018),  427--450.

\bibitem{Hac}  Ian Hacking, 
\emph{Why Is There Philosophy of Mathematics at All?},
Cambridge University Press, 2005.

\bibitem{Hal} T.C. Hales, 
\emph{A proof of the Kepler conjecture},
 Annals of Mathematics, 162 (2005), 1065--1185.
 
 \bibitem{Hal3}   T. Hales,
\emph{Formal proof},
 Notices Amer Math Soc 55 (2008), 1370--1381.
 

 
 \bibitem{Hal2} T.C. Hales and 21 others, 
\emph{A formal proof of the Kepler conjecture},
 Forum of Mathematics, Pi 5, e2 (2017), 29 pages.
 
 \bibitem{Hall} Michael Hallett, 
\emph{Foundations of Mathematics},
The Cambridge History of Philosophy, 1870--1945, ed. Thomas Baldwin
 (Ch.~10 of Part I: 1870--1914), 
Cambridge University Press (2003), 128--156, 833--837.


 \bibitem{Hara}
Donna Haraway,
\emph{Situated knowledges: The science question in feminism and the privilege of partial perspective}, 
Feminist studies, 14 (1988), 575--599.

 \bibitem{Harris} Michael Harris's blog:
  \emph{Mathematics without apologies}, 
  June 2nd, 2018, {\tt https://mathematicswithoutapologies.wordpress.com}



\bibitem{Harr}   John Harrison,
\emph{Formal proof -- theory and practice},
 Notices Amer Math Soc 55 (2008), 1395--1406.

 

 \bibitem{Hosa}
John Hosack, 
\emph{An inclusive philosophy of mathematics},
Notices Amer. Math. Soc. 66 (2019),  1433--1437.

 


\bibitem{Kah}
Reinhard Kahle
\emph{Towards the Structure of Mathematical Proof},
CEUR Workshop proceedings vol 1186 (2014), paper \# 22.

\bibitem{Kah2}
Reinhard Kahle
\emph{What is a proof?}
Axiomathes  25 (2015), 79--91.

\bibitem{Knu}
Donald Knuth
\emph{All Questions Answered},
Notices of the AMS,  49  (2002), 318--324.

 \bibitem{Laba} 
 Izabella Laba
\emph{Photography},
in ``Art in the life of mathematicians'' (ed. Anna Szemer\'edi), 
American Mathematical Society, 2015.

 \bibitem{Lag} 
 J. Lagarias, 
\emph{The Kepler Conjecture and its Proof},
 (Springer, 2009), 3--26.
 
 
  \bibitem{Lang}
M. Lange,
\emph{Because Without Cause: Non-causal Explanations in Science and Mathematics}, 
2017, Oxford: Oxford University Press.

\bibitem{Ler} 
X. Leroy, 
\emph{Formal certification of a compiler back-end,
or: programming a compiler with a proof assistant}, 
33rd ACM Symposium on \emph{Principles of Programming Languages}, 
ACM Press (2006)  42--54.
   

 
  \bibitem{Man1}
 Paolo Mancosu,
\emph{Mathematical Explanation: Problems and Prospects},
Topoi 20 (2001),  97--117.

 \bibitem{Man2}
Paolo Mancosu,
\emph{Explanation in Mathematics}, The Stanford Encyclopaedia of Philosophy (2018), Edward N. Zalta (ed.),  https://plato.stanford.edu/archives/sum2018/entries/mathematics-explanation/

 \bibitem{Mass}
Patrick Massot,
\emph{Formal mathematics for mathematicians and mathematics students}
lecture on youtube.com (https://www.youtube.com/watch?app=desktop\&v=tp\underline{ }h3vzkObo).

 \bibitem{May}
James Maynard,
\emph{Small gaps between primes},
Ann. of Math.   181 (2015),  383--413.

\bibitem{Mit}  Melanie Mitchell,
\emph{Artificial Intelligence:  A Guide for Thinking Humans},
Farrar, Straus and Giroux, 2019.


 \bibitem{Mord1}  L. J. Mordell,
 \emph{On the rational solutions of the indeterminate equations of the third and fourth degrees},
 Proc. Cab. Phil. Soc 21 (1922), 179--192.

 \bibitem{Mord2}  L. J. Mordell,
 \emph{Reminiscences of an Octogenarian Mathematician},
The American Mathematical Monthly  78 (1971),  952--961

\bibitem{Nath08}
Melvyn B. Nathanson,   
\emph{Desperately seeking mathematical truth},
Notices of the AMS (2008),  773.

\bibitem{Nath09}
Melvyn B. Nathanson,   
\emph{Desperately seeking mathematical proof},
Math. Intelligencer 31 (2009),   8--10.

\bibitem{PaHa77}
J. Paris, and L. Harrington,
\emph{A Mathematical Incompleteness in Peano Arithmetic}
 In Barwise, J. (ed.). Handbook of Mathematical Logic. Amsterdam, Netherlands: North-Holland, 1977.
 
\bibitem{RTV}
 Colin J. Rittberg, Fenner S. Tanswell and Jean-Paul Van Bendegem,
\emph{Epistemic injustice in mathematics}
Synthese 197 (2020), 3875--3904


\bibitem{Rob15}
Siobhan Roberts,
\emph{In Mathematics, Mistakes Aren't What They Used To Be},
Nautilus online, 2015.

\bibitem{Rock}
Dan Rockmore,
\emph{Prove It!},
The New York Review of Books, Jan 13, 2022.

\bibitem{RSST} Neil Robertson, Daniel P. Sanders, Paul  Seymour, and  Robin Thomas,  
\emph{A new proof of the four-colour theorem},
Electron. Res. Announc. Amer. Math. Soc. 2 (1996),  17--25.

\bibitem{Slim} 
Dirk Schlimm,  
\emph{Peano on Symbolization, Design Principles for Notations, and the Dot Notation},
Philosophia Scientae, 25 (2021), 95--126.

\bibitem{Sch}
Peter Scholze,
 \emph{Half a year of the Liquid Tensor Experiment: Amazing developments},
https://xenaproject.wordpreWhy abc is still a conjecture.com/2021/06/05 

\bibitem{Sch2}
Peter Scholze (with Dustin Clausen),
 \emph{Lectures on Analytic Geometry},
https://www.math.uni-bonn.de/people/scholze/Analytic.pdf

\bibitem{SS}
Peter Scholze and Jakob Stix,
\emph{Why $abc$ is still a conjecture} (preprint, 2018).

\bibitem{Tarski} 
Alfred Tarski,
\emph{Truth and Proof},
Scientific American 220  (1969), 63--77.


\bibitem{Th84} 
William P.  Thurston,
\emph{On proof and progress in mathematics},
 Bull. Amer. Math. Soc. (N.S.) 30 (1994), no. 2, 161--177
 


\bibitem{Tur1} A. M. Turing,
\emph{Computing Machinery and Intelligence},
Mind, New Series,   59  (1950),  433--460.

\bibitem{Venk}
Akshay Venkatesh,
\emph{How we place value in mathematics}, preprint.


\bibitem{Voev}
Vladimir Voevodsk\u{\i},
\emph{What If Current Foundations of Mathematics Are Inconsistent?}
Lecture at https://www.ias.edu/ideas/2012/voevodsky-foundations-of-mathematics, 2010.


 


\bibitem{Wied}   Freek Wiedijk,
\emph{Formal proof -- getting started},
 Notices Amer Math Soc 55 (2008), 1408--1414.

   
\end{thebibliography}

\end{document}